\documentclass[a4paper,11pt,reqno,twoside]{amsart}

\usepackage{amsmath}
\usepackage{amsfonts}
\usepackage{amssymb}
\usepackage{amsthm}
\usepackage{color}
\usepackage{ifpdf}
\usepackage{array}
\usepackage{url}
\usepackage{multirow,bigdelim}
\usepackage{ascmac}
\usepackage{mathrsfs}
\usepackage{ulem}

\newtheorem{theorem}{Theorem}[section]
\newtheorem{proposition}[theorem]{Proposition}
\newtheorem{lemma}[theorem]{Lemma}
\newtheorem{corollary}[theorem]{Corollary}

\newtheorem{remark}[theorem]{Remark}
\newtheorem{example}[theorem]{Example}

\numberwithin{equation}{section}

\newcommand*{\C}{\mathbb{C}}
\newcommand*{\R}{\mathbb{R}}
\newcommand*{\Z}{\mathbb{Z}}
\newcommand{\N}{\mathbb N}
\newcommand{\D}{\mathscr D}
\newcommand{\Dom}{\mathfrak D}
\newcommand{\ran}{\operatorname{ran}}

\newcommand{\PW}{\mathrm{PW}}
\newcommand{\Hglob}{\mathcal H_W}
\newcommand{\Ha}{\mathcal H(A_a)}
\newcommand{\Hhat}{\widehat{\mathcal H(A_a)}}
\newcommand{\Hstar}{\mathcal H_\ast}
\newcommand{\Samp}{\mathcal S}
\newcommand{\Jref}{\mathcal R}
\newcommand{\Dlim}{\mathsf D_{\pi/2}}
\newcommand{\Wdist}{\mathscr W_C}

\newcommand{\Hred}{\mathcal H^{\circ}}

\newcommand{\Psf}{\Delta_\Gamma}

\newcommand{\comment}[1]{}
\title[A spectral theory of zeta functions of function fields]{
A spectral theory of \\
zeta functions of function fields
}
\author[M. Suzuki]{Masatoshi Suzuki}
\date{Version of \today}
\AtBeginDocument{%
\begin{abstract}
We study Hilbert spaces constructed from 
the Weil quadratic form associated with the zeta function of a curve over a finite field, 
together with self-adjoint extensions of the differentiation operator on these spaces. 
For a distinguished family of these extensions, 
we prove that the eigenvalues converge to 
the ordinates of the zeros of the zeta function and that the corresponding entire spectral functions, 
after normalization, converge locally uniformly to an infinite product determined by these zeros.
These limits occur as the interval parameter approaches a finite endpoint, 
rather than infinity as in the setting of the Riemann zeta function.
We also realize the global Hilbert space through 
the Rosati trace pairing and show that translation 
by $\log q$ is unitarily equivalent to multiplication by normalized Frobenius.
This provides a geometric model for 
the analogous analytic constructions associated with the Riemann zeta function.
\end{abstract}
\maketitle
}
\begin{document}

\section{Introduction} \label{section_1}

The nontrivial zeros of the Riemann zeta function are fundamental objects
of study in analytic number theory, particularly in the theory of the
distribution of prime numbers.
In \cite{Su26b}, assuming the Riemann Hypothesis (RH), we proved that
their ordinates can be approximated by eigenvalues of self-adjoint
extensions of a differentiation operator on Hilbert spaces obtained by
completing spaces of functions on finite intervals with respect to the
local Weil quadratic form.
That result is based on the study of the local Weil quadratic form in \cite{Su26a}.
\medskip

Weil's explicit formula generalizes the classical explicit formulas of
analytic number theory.
The Weil distribution appearing in this formula gives rise to the Weil
quadratic form.
Since its introduction in \cite{We52}, the Weil distribution has provided
a unified framework for treating zeta functions of both number fields and
function fields.
It is therefore natural to consider a function field analogue of the theory
developed in \cite{Su26a,Su26b}.
Moreover, since RH has been proved for zeta functions of function fields
\cite[Section~1]{Milne}, one may expect unconditional analogues of results
that were obtained under RH in the case of the Riemann zeta function.
Although the basic methods are similar, some of the results are analogous 
in the two settings, whereas others are specific to function fields.

The study of localized Weil quadratic forms, 
originating in the variational work of Yoshida and Bombieri, 
has developed into an operator-theoretic approach involving associated 
Hilbert spaces and self-adjoint operators \cite{CC23,CCM25,Su26a}. 
Another direction concerns the relation of Weil positivity to geometry and 
trace formulas \cite{CC21,CCM09}. 
We follow the former approach in the function field setting. 
The present paper is motivated by \cite{Su26a,Su26b}, 
but is self-contained and can be read independently of those papers.

More precisely, we show that the ordinates of the nontrivial zeros of a
function field zeta function can be approximated by eigenvalues of suitable
self-adjoint extensions of a differentiation operator.
Let $C/\mathbb F_q$ be a geometrically connected nonsingular projective
curve of genus $g\ge1$.
Writing $N_n=\#C(\mathbb F_{q^n})$, the zeta function is defined by
\[
 Z_C(T)=\exp\left(\sum_{n\ge1}\frac{N_n}{n}T^n\right)
       =\frac{P_C(T)}{(1-T)(1-qT)}.
\]
The numerator \(P_C\) can be determined from \(N_1,\ldots,N_g\). 
Thus its zeros can, in principle, be computed by a finite algebraic procedure.
The purpose of the present work is not to provide another method for
computing those zeros, but to construct the local spectral theory
arising from the Weil quadratic form and to identify the action of
Frobenius within the resulting operator theory.
\medskip

We define
\begin{equation} \label{eq_101}
 \xi_C(s):=q^{g(s-1/2)}P_C(q^{-s}),\qquad
 A_C(z):=\xi_C(1/2-iz).
\end{equation}
By Weil's theorem for the zeta function of a curve, $A_C$ is a real even
entire function, and all its zeros are real \cite[Section~1]{Milne}.

Let $\Gamma$ be the set of distinct zeros of $A_C$, and put $h=\log q$.
Choose representatives $\beta_1,\ldots,\beta_d$ of the $d$ distinct
residue classes of zeros modulo $2\pi/h$.
Then
\[
 \Gamma=\bigsqcup_{j=1}^d\left(\beta_j+\frac{2\pi}{h}\Z\right).
\]
Set $\Gamma_+:=\Gamma\cap(0,\infty)$, 
\begin{equation} \label{eq_102}
 \chi:=\begin{cases}0,&0\notin\Gamma,\\1,&0\in\Gamma,\end{cases}
 \qquad \varepsilon:=(-1)^\chi,
\end{equation}
and define
\begin{equation} \label{eq_103}
 \Psf(z):=z^\chi\prod_{\gamma\in\Gamma_+}
              \left(1-\frac{z^2}{\gamma^2}\right).
\end{equation}
Periodicity implies that
$\sum_{\gamma\in\Gamma_+}\gamma^{-2}<\infty$, so the product converges
locally uniformly on the whole plane.
The condition $0\in\Gamma$ is equivalent to $P_C(q^{-1/2})=0$;
we refer to this as the case in which there is a central zero.
The function $\Psf$ has a simple zero at each point of $\Gamma$ and does
not record the multiplicities of the zeros of $A_C$.
We define the endpoint by
\[
 a_\ast:=\frac{dh}{2},
\]
which depends on the number $d$ of distinct residue classes of zeros, 
rather than on the genus alone.

We define the Fourier transform by 
$\widehat u(z):=\int_\R u(x)e^{izx}\,dx$. 
Let $m_\gamma$ denote the multiplicity of $\gamma\in\Gamma$ as a zero of $A_C$.
The Weil quadratic form on $C_c^\infty(\R)$ is defined by
\begin{equation} \label{eq_104}
 Q_C(u,f):=\sum_{\gamma\in\Gamma}m_\gamma\widehat u(\gamma)
                  \overline{\widehat f(\gamma)},\qquad
 Q_C(u):=Q_C(u,u).
\end{equation}
Since all zeros of $A_C$ are real, $Q_C$ is nonnegative.
By Weil's explicit formula, the Weil quadratic form and 
its representing operator on any fixed bounded interval 
can be constructed from finitely many values of $N_n$, 
without computing the zeros of $P_C$ (see Section~\ref{section_3}).
Let $A_a$ be the nonnegative self-adjoint operator on $L^2(-a,a)$ representing
the bounded extension of $Q_C|_{C_c^\infty(-a,a)}$.
We denote by $\Ha$ the completion of $C_c^\infty(-a,a)$ with respect to $Q_C$.
For $0<a\le a_\ast$, the $Q_C$-norm is equivalent to the usual
$L^2$-norm, and $\Ha$ coincides with $L^2(-a,a)$ as a set.
On the other hand, for $a>a_\ast$, $Q_C$ is degenerate on $L^2(-a,a)$.
In particular, $Q_C$ remains nondegenerate at the endpoint.
These facts will be proved in Section~\ref{section_4}.

For $0<a\le a_\ast$, define the differentiation operator
$\D_a=i\,d/dx$ on the Hilbert space $\Ha$, with domain
$C_c^\infty(-a,a)$.
It is symmetric with deficiency indices $(1,1)$ and has a family of
self-adjoint extensions $\overline\D_{a,\theta}$ parametrized by a phase
$\theta\in[0,2\pi)$; see Section~\ref{section_5}.
The following phase plays an important role throughout this paper:
\begin{equation} \label{eq_105}
\theta_\ast
=
\begin{cases}
\pi,&0\not\in\Gamma,\\ 0,&0\in\Gamma
\end{cases}.
\end{equation}
In particular, $e^{i\theta_\ast}=-\varepsilon$.

\begin{theorem} \label{thm_1}
Write
$\Gamma_+=\Gamma\cap(0,\infty)=\{\gamma_1<\gamma_2<\cdots\}$.
For $0<a\le a_\ast$, let $0<t_{a,1}<t_{a,2}<\cdots$ be the positive
eigenvalues of the self-adjoint operator
$\overline\D_{a,\theta_\ast}$.
Then, for each $n\ge1$,
\begin{equation}
 \label{eq_106}
 t_{a,n}\ge\gamma_n,\qquad
 t_{a,n}\searrow\gamma_n\quad(a\uparrow a_\ast),\qquad
 t_{a_\ast,n}=\gamma_n.
\end{equation}
Here $\searrow$ denotes monotone nonincreasing convergence.
At $a=a_\ast$, the spectrum of $\overline\D_{a,\theta_\ast}$ coincides
with $\Gamma$, and every eigenvalue is simple.
If there is a central zero, zero is a simple eigenvalue for every $a$.
\end{theorem}

The multiplicities of the zeros of the zeta function must be distinguished
from those of the eigenvalues of the differentiation operator.
Simplicity of the latter does not require simplicity of the former.

For each $z\in\C$, the Riesz vector $u_{z,a}$ representing the continuous
linear functional $u\mapsto\widehat u(\bar z)$ on $\Ha$ is determined by
\[
 \widehat u(\bar z)=\langle u,u_{z,a}\rangle_{A_a}
 \quad(:=Q_C(u,u_{z,a})).
\]
Using the vectors $u_{i,a}$ and $u_{-i,a}$, which belong to the
deficiency subspaces of $\D_a$ and have equal norms, we define the
spectral function associated with the boundary form by
\[
 W(a,\theta;z):=(z-i)\widehat{u_{i,a}}(z)
       +e^{i\theta}(z+i)\widehat{u_{-i,a}}(z).
\]
This is an entire function whose zeros coincide with the eigenvalues of
the self-adjoint extension $\overline\D_{a,\theta}$
(Theorem~\ref{thm_502}).

\begin{theorem} \label{thm_2}
We have
\begin{equation} \label{eq_107}
 \lim_{a\uparrow a_\ast}\frac{W(a,\theta_\ast;z)}{W(a,\theta_\ast;i)}
       =\frac{\Delta_\Gamma(z)}{\Delta_\Gamma(i)},
\end{equation}
where the convergence is uniform on every compact subset of $\C$.
In particular, if all zeros of $A_C$ are simple, then $\theta_\ast=\pi$
and the right-hand side is $A_C(z)/A_C(i)$.
\end{theorem}

If $A_C$ has a central zero, then $W(a,\theta_\ast;0)=0$.
Thus, unlike the normalization at the origin in \cite{Su26b}, we
normalize at $i$ in \eqref{eq_107}, 
since this choice works in both cases.
The limit of $W(a,\theta_\ast;z)$ itself can also be computed explicitly.
Furthermore, after multiplication by a holomorphic factor independent of $a$,
we obtain a limit formula for $A_C/(A_C+iA_C')$ in the upper half-plane
(Theorem~\ref{thm_607} and Corollary~\ref{cor_608}).

Sections~\ref{section_8} and \ref{section_9} provide
a geometric and operator-theoretic interpretation of the
self-adjoint differentiation operator at $a=a_\ast$
appearing in Theorem~\ref{thm_1}.
More importantly, they give a geometric interpretation, in the function field
setting, of part of the analytic construction developed in
\cite{Su26a,Su26b}.
The global Hilbert space is realized through the Rosati trace pairing,
and translation by $\log q$ is unitarily equivalent to multiplication
by normalized Frobenius.
\medskip

For ease of comparison, the order of the main theorems follows that of
\cite{Su26b}, which treats the Riemann zeta function. 
In the function field case, however, Theorem~\ref{thm_2} is proved independently of Theorem~\ref{thm_1}.
This difference reflects the different structures of the families of
Hilbert spaces defined by the Weil quadratic form in the two settings.
\medskip

The paper \cite{Su26a}, on which \cite{Su26b} is based, studies the
self-adjoint operator $A_a$ representing the Weil quadratic form
associated with the Riemann zeta function on $L^2(-a,a)$.
Its results include a realization as a Friedrichs extension, continuity
of the lowest eigenvalue $\lambda_a$ for $a>0$, an asymptotic formula
and simplicity for sufficiently small $a>0$, and evenness of the
corresponding eigenfunctions
\cite[Theorems~1.1, 1.3, 1.4 and Corollary~1.2]{Su26a}.
We also consider the corresponding questions for function fields.

In the function field case, $A_a$ is bounded, and the Friedrichs
extension of the corresponding symmetric operator is simply its closure.
The lowest eigenvalue of $A_a$, however, is generally discontinuous and
varies as a step function.
Moreover, $A_a=2ghI$ for $0<2a\le h$, so the eigenspace corresponding
to the lowest eigenvalue is infinite-dimensional.
These spectral phenomena differ from those in \cite{Su26a}; see
Sections~\ref{section_3}--\ref{section_4}.
\medskip

Finally, we discuss the historical context of Weil positivity and its
relation to geometric positivity in the function field setting.
Weil's geometric proof of RH for function fields must be distinguished
from his positivity criterion based on the explicit formula of 1952
\cite{We52}.
The former uses intersection theory for correspondences on the surface
$C\times C$ \cite[Section~1]{Milne}.
The latter formulates positivity of the Weil distribution as an analytic
condition equivalent to RH, applicable to both number fields and
function fields.
In fact, in \cite{We52}, Weil deduces positivity of the corresponding
distribution from the already established RH for function fields;
he does not derive that positivity directly from intersection theory.

Concerning the relation between geometric positivity and Weil positivity,
Swinnerton-Dyer~\cite[Section~7]{SD71} constructs a map from
correspondences on a curve to endomorphisms of its Jacobian and proves
positivity of the trace pairing with respect to transposition by means of
intersection theory.
Patterson~\cite[Section~5.17]{Pa88} explains an argument for geometric
positivity using Frobenius on the Jacobian and the Rosati involution.
Bombieri~\cite[Section~5]{Bom06} also states that Weil positivity for
function fields follows from the algebraic index theorem for surfaces,
but does not give an explicit derivation there.
Tretkoff~\cite[Section~2.2]{Tr06} formulates the discrete Weil quadratic
form for function fields and explains its relation to positivity of the
Rosati involution.
Connes, Consani, and Marcolli~\cite[Section~2.3]{CCM09} discuss
positivity arising from intersection theory for correspondences;
in \cite[Sections~6.2--7]{CCM09}, they examine the equivalence between
RH and positivity of a trace pairing, together with its geometric
interpretation.
The monograph of van Frankenhuijsen~\cite[Chapters~4 and~6]{vF14}
treats the analytic formulation of Weil positivity and an approach
through traces of shift operators, but the proof of positivity by the
latter approach remains incomplete.

These references do not make explicit the precise relation between 
the quadratic form $Q_C$ considered here and positivity arising from intersection
theory.
In Section~\ref{section_8}, we therefore express $Q_C$ as an
integral of the Rosati trace pairing, following the arguments of
Swinnerton-Dyer, Patterson, and Tretkoff, and explain its relation to
geometric positivity.
\medskip

The paper is organized as follows.
Section~\ref{section_2} introduces abstract periodic zero data and 
explains how they arise from zeta functions of curves over finite fields.
In Section~\ref{section_3}, we use Weil's explicit formula to construct the associated screw function and the bounded self-adjoint operators representing the local Weil quadratic form.
Section~\ref{section_4} studies the positivity of the local Weil quadratic form, the spectrum of its representing operators, and the sampling map at the endpoint.
Section~\ref{section_5} constructs self-adjoint extensions of the differentiation operator on the resulting Hilbert spaces and introduces their entire spectral functions.

In Section~\ref{section_6}, we determine the reproducing kernel at the endpoint and prove the main theorems concerning the convergence of eigenvalues and spectral functions.
We also obtain limit formulas involving the zeta function of the curve.
Section~\ref{section_7} constructs the global Hilbert space, describes its realizations in terms of de Branges spaces and model spaces, and studies the convergence of compressed inverse operators and resolvents.

Section~\ref{section_8} expresses the Weil quadratic form in terms of the Rosati trace pairing on the endomorphism algebra of the Jacobian and explains its relation to geometric positivity.
Using this realization, Section~\ref{section_9} identifies the global differentiation operator as the self-adjoint generator of a unitary translation group and shows how normalized Frobenius appears in this representation.
Section~\ref{section_10} compares the resulting theory with that developed for the Riemann zeta function.
Finally, Appendix~\ref{section_app} presents explicit computations for several curves over finite fields.

\section{Abstract periodic data and zeta functions of curves} \label{section_2}

Throughout this paper, inner products are linear in the first variable.
We write $F^\sharp(z):=\overline{F(\bar z)}$ and
$\C_+:=\{z:\Im z>0\}$.
The Fourier transform convention is
\begin{equation} \label{eq_201}
 (\mathcal Fu)(z)=\widehat u(z):=\int_\R u(x)e^{izx}\,dx,
 \qquad
 \|\mathcal Fu\|_{L^2}^2=2\pi\|u\|_{L^2}^2.
\end{equation}
Functions on an interval are regarded as functions on $\R$ by extension
by zero.
We put
\[
 \PW_a:=\mathcal F L^2(-a,a).
\]
In particular, $\widehat{iu'}(z)=z\widehat u(z)$ for
$u\in C_c^\infty(\R)$.

\subsection{Abstract periodic data}

We consider abstract periodic data
$C:=(\Gamma,(m_\gamma)_{\gamma\in\Gamma})$ of the following form,
which include the zero sets of zeta functions of curves as special cases.
We also use the notation $C$ for these data to allow a direct comparison with the case of a curve. 
Choose $h>0$, real numbers $\beta_j$ distinct modulo $2\pi/h$, and
positive integers $m_j$ $(1\le j\le d)$, and set
\begin{equation} \label{eq_202}
 \Gamma:=\bigsqcup_{j=1}^d\left(\beta_j+\frac{2\pi}{h}\Z\right),\qquad
 m_{\beta_j+2\pi n/h}:=m_j.
\end{equation}
Assume that $\Gamma=-\Gamma$, $m_{-\gamma}=m_\gamma$, and that,
for a nonzero real constant $c_A$, the function
\begin{equation} \label{eq_203}
 A(z):=c_A\prod_{j=1}^d\sin^{m_j}\frac{h(z-\beta_j)}2
\end{equation}
is even.
For these abstract data, we again define $\Psf$ by the product over
pairs of positive and negative zeros in \eqref{eq_103}, and let $\chi$
be the indicator of $0\in\Gamma$, as in \eqref{eq_102}.
Writing
\[
 S(z):=\prod_{j=1}^d\sin\frac{h(z-\beta_j)}2,
\]
we obtain the following sine product representation of $\Psf$, together
with the exponential types of $\Psf$ and $A$:
\begin{equation} \label{eq_204}
 \begin{gathered}
 \Psf(z)=\begin{cases}S(z)/S(0),&0\notin\Gamma,\\
                         S(z)/S'(0),&0\in\Gamma,\end{cases}\\
 a_\ast:=\frac{dh}{2},\qquad
 \tau_A:=\frac h2\sum_{j=1}^d m_j.
 \end{gathered}
\end{equation}
Indeed, $S$ has exponential type $a_\ast$ and has precisely the zeros
$\Gamma$, all of which are simple.
By Hadamard factorization and parity, $S$ is a nonzero real constant multiple of $\Psf$. 
Its value or derivative at the origin determines the constant.
Shifting a representative $\beta_j$ by an integer multiple of $2\pi/h$
changes $S$ only by a sign, so \eqref{eq_204} is independent of
the choice of representatives.
By the definition of $\Psf$, we have
\[
 \Psf(-z)=\varepsilon\Psf(z).
\]
Since $A$ is even, its zero at the origin has even multiplicity if it exists.
We write $m_0$ for this multiplicity if $0\in\Gamma$, and set $m_0=0$
otherwise.

\subsection{Specialization to the zeta function of a curve}

For a curve, rationality, the functional equation, and Weil's theorem give
\[
 P_C(T)=\prod_{\nu=1}^{2g}(1-\omega_\nu T)\in\Z[T],\quad
 P_C(T)=q^gT^{2g}P_C(1/(qT)),\quad
 |\omega_\nu|=\sqrt q;
\]
see \cite[Section~1]{Milne}.
The functional equation implies $\xi_C(s)=\xi_C(1-s)$, and the real
coefficients of $P_C$ imply that $A_C$ is a real even entire function.
As before, put $h=\log q$.
For the distinct phases in $\omega_\nu=\sqrt q e^{i\theta_\nu}$,
choose $\beta_j\equiv-\theta_\nu/h\pmod{2\pi/h}$.
This gives \eqref{eq_202} and \eqref{eq_203}.
Indeed, applying
\[
 1-e^{ih(z-\beta_j)}=-2i\,e^{ih(z-\beta_j)/2}
                 \sin\frac{h(z-\beta_j)}2
\]
to each factor, the exponential factors depending on $z$ cancel the
prefactor in \eqref{eq_101}.
The remaining constant is nonzero and real because $A_C$ is real on the
real axis, and $\sum_jm_j=2g$.
In the case of a curve, we henceforth take $A=A_C$.
\medskip

The function $A_C$ is a nonzero trigonometric polynomial and does not
belong to $L^2(\R)$.
Thus, although its exponential type is $\tau_A=gh$, it does not belong
to $\PW_{\tau_A}$.
If there are multiple zeros, then $a_\ast<\tau_A$.
We exclude genus zero, since the measure arising from the zeros and the
Hilbert spaces below are then trivial, but we do not exclude central or
multiple zeros when $g\ge1$.

\section{Representation of the Weil quadratic form} \label{section_3}

This section treats the function field analogues of
\cite[Theorem 1.1, Corollary 1.2]{Su26a}.

\subsection{Weil's explicit formula for function fields}

For the abstract periodic data $C$ introduced in the preceding section
and $u,f\in C_c^\infty(\R)$, define $Q_C(u,f)$ by \eqref{eq_104}.
Periodicity gives $\sum_\gamma m_\gamma/(1+\gamma^2)<\infty$,
so the rapid decay of the Fourier transforms implies absolute convergence.
Put
\[
 p_n:=\sum_{j=1}^{d}m_j\zeta_j^n,\qquad
 \zeta_j:=e^{ih\beta_j}.
\]
Then $p_{-n}=p_n\in\R$.
The Poisson summation formula gives
\begin{equation} \label{eq_301}
 \Wdist(\phi):=\sum_{\gamma\in\Gamma}m_\gamma\widehat\phi(\gamma)
       =h\sum_{n\in\Z}p_n\phi(nh),\qquad
 \phi\in C_c^\infty(\R).
\end{equation}
Indeed, multiply
$\sum_{k\in\Z}e^{2\pi ikx/h}=h\sum_{n\in\Z}\delta_{nh}$
by $m_je^{i\beta_jx}$ and sum over $j$.
For each test function, the right-hand side is a finite sum.
In the case of a curve, $p_0=2g$ and
\begin{equation} \label{eq_302}
 p_n=\frac{q^n+1-N_n}{q^{n/2}},\qquad
 N_n=\#C(\mathbb F_{q^n}),\qquad(n\ge1).
\end{equation}
Thus \eqref{eq_301} is Weil's explicit formula in terms of the
numbers of rational points of the curve.
Writing $\widetilde f(x)=\overline{f(-x)}$, we have
\[
 Q_C(u,f)=\Wdist(u\ast\widetilde f).
\]

\subsection{The screw function for function fields}

A continuous function $g$ on $\R$ satisfying $g(t)=\overline{g(-t)}$ is called a
screw function on the real line in the sense of Krein and
Langer~\cite[Section 5]{KrLa14} if the kernel
$g(t-u)-g(t)-g(-u)+g(0)$ is nonnegative definite on $\R$.

\begin{proposition} \label{prop_301}
The real even function defined by 
\begin{equation} \label{eq_303}
 g_C(t):=-\frac{m_0}{2}t^2+
       \sum_{\substack{\gamma\in\Gamma\\\gamma\ne0}}
       m_\gamma\frac{\cos(\gamma t)-1}{\gamma^2}
\end{equation}
is a continuous screw function and admits the representation
\begin{equation} \label{eq_304}
 g_C(t)=-\frac{hp_0}{2}|t|
        -h\sum_{n\ge1}p_n(|t|-nh)_+,
\end{equation}
where $(x)_+=\max(x,0)$, and the sum on the right is locally finite.
Moreover,
\begin{equation} \label{eq_305}
 -g_C''=h\sum_{n\in\Z}p_n\delta_{nh}
\end{equation}
in the sense of distributions.
\end{proposition}

\begin{proof}
Since $\sum_\gamma m_\gamma/(1+\gamma^2)<\infty$, the series in
\eqref{eq_303} converges uniformly on $\R$.
By the reflection symmetry of $\Gamma$, replacing the cosine by
$e^{-i\gamma t}$ does not change the sum.
Hence
\[
 g_C(t-s)-g_C(t)-g_C(-s)+g_C(0)
 =m_0ts+\sum_{\gamma\ne0}\frac{m_\gamma}{\gamma^2}
       (e^{-i\gamma t}-1)(e^{i\gamma s}-1).
\]
The right-hand side is a sum of nonnegative definite kernels, so $g_C$
satisfies the defining condition for a screw function.
Differentiating \eqref{eq_303} twice in the sense of
distributions and using the reflection symmetry of $\Gamma$, we obtain
\[
 -g_C''=\sum_{\gamma\in\Gamma}m_\gamma e^{i\gamma t}.
\]
Interpreting the right-hand side as a distribution,
\eqref{eq_301} yields \eqref{eq_305}.
On the other hand,
\[
 \frac{d^2}{dt^2}|t|=2\delta_0,\qquad
 \frac{d^2}{dt^2}(|t|-a)_+
   =\delta_a+\delta_{-a}\quad(a>0)
\]
in the sense of distributions.
If $G(t)$ denotes the right-hand side of \eqref{eq_304},
then
\[
 -G''=hp_0\delta_0+
 h\sum_{n\ge1}p_n(\delta_{nh}+\delta_{-nh})
 =h\sum_{n\in\Z}p_n\delta_{nh}.
\]
Thus $(g_C-G)''=0$.
Since both $g_C$ and $G$ are even and vanish at the origin, $g_C=G$.
\end{proof}

For a curve, \eqref{eq_302} and \eqref{eq_304}
determine $g_C$ on each bounded interval from the numbers of rational
points.
In particular, $g_C(t)=-gh|t|$ for $|t|\le h$.
Unlike the $|t|\log|t|$ singularity in the Riemann zeta case, the term at
the origin is piecewise linear.

\subsection{Local operators and the Friedrichs extension}

Applying \eqref{eq_301} to $u\ast\widetilde f$, we obtain, for each
$a>0$, an extension of $Q_C$ to a bounded nonnegative quadratic form on
$L^2(-a,a)$.
The corresponding bounded nonnegative self-adjoint operator,
characterized by
\[
 Q_C(u,f)=\langle A_a u,f\rangle_{L^2},\quad
 u,f\in L^2(-a,a),
\]
is given by
\begin{equation} \label{eq_306}
 (A_a u)(x)=h\left(p_0u(x)+\sum_{1\le n<2a/h}p_n
       \{u(x+nh)+u(x-nh)\}\right),\qquad |x|<a.
\end{equation}
We use extension by zero and omit shifts for which the intervals meet
only at an endpoint. Consequently,
\begin{equation} \label{eq_307}
 0\le Q_C(u)\le
 h\left(p_0+2\sum_{1\le n<2a/h}|p_n|\right)\|u\|_{L^2}^2,
 \qquad u\in L^2(-a,a).
\end{equation}
Let $K_a$ be the integral operator on $L^2(-a,a)$ with kernel
$g_C(x-y)$, and let $P_a$ be the usual $L^2$-orthogonal projection onto
\[
 L_0^2(-a,a)=
 \{f\in L^2(-a,a)~\mid~\int_{-a}^{a}f(x)dx=0\}.
\]
Put $G_a=P_aK_aP_a$.
On the usual $L^2$ space, set
\[
 d_a=i\,d/dx,\quad \Dom(d_a)=H_0^1(-a,a),\qquad
 d_a^\ast=i\,d/dx,\quad\Dom(d_a^\ast)=H^1(-a,a),
\]
where $\Dom(T)$ denotes the domain of an operator $T$.

\begin{proposition} \label{prop_302}
Define $B_a:=d_a^\ast G_a d_a$ on $L^2(-a,a)$ with
$\Dom(B_a)=H_0^1(-a,a)$. Then
\begin{equation} \label{eq_308}
 B_a=A_a|_{H_0^1(-a,a)},\qquad \overline{B_a}=A_a.
\end{equation}
Thus $B_a$ is essentially self-adjoint, and its Friedrichs extension is $A_a$.
The space $C_c^\infty(-a,a)$ is a core for the quadratic form associated
with $A_a$, and
\begin{equation} \label{eq_309}
 \inf\sigma(A_a)=
 \inf_{0\ne u\in C_c^\infty(-a,a)}
 \frac{Q_C(u)}{\|u\|_{L^2}^2}.
\end{equation}
\end{proposition}
\begin{proof}
Since $g_C$ is Lipschitz continuous on every compact interval,
$K_a:L^2(-a,a)\to H^1(-a,a)$ is bounded.
For $u\in H_0^1(-a,a)$, the function $d_a u$ has mean zero, and
integration by parts gives $K_a d_a u=i(g_C'\ast u)$.
The constant subtracted by $P_a$ disappears upon differentiation, so
\eqref{eq_305} yields
\[
 d_a^\ast G_a d_a u=-g_C''\ast u=A_a u.
\]
Thus $B_a$ is well defined on the stated domain $H_0^1(-a,a)$.
In particular,
\[
 Q_C(u)=\langle G_a d_a u,d_a u\rangle_2
       =\int_{-a}^a\int_{-a}^a
                 g_C(x-y)u'(y)\overline{u'(x)}\,dy\,dx.
\]
Since $A_a$ is bounded and $H_0^1(-a,a)$ is dense in $L^2(-a,a)$,
the closure of its restriction $B_a$ is $A_a$.
The correspondence between closed nonnegative forms and self-adjoint
operators shows that its Friedrichs extension is the same operator.
Moreover, $(Q_C(u)+\|u\|_2^2)^{1/2}$ is equivalent to the $L^2$-norm.
The usual approximation by test functions therefore shows that
$C_c^\infty(-a,a)$ is a core for the quadratic form and gives
\eqref{eq_309}.
\end{proof}

\section{Lower bounds for the local Weil quadratic form} \label{section_4}

We investigate the function field counterparts of the questions addressed
in \cite[Theorems~1.3 and 1.4]{Su26a}.

\subsection{Completion with respect to the Weil quadratic form}

Put $N_a:=\{u\in C_c^\infty(-a,a):Q_C(u)=0\}$ and define
\[
 \Ha:=\overline{C_c^\infty(-a,a)/N_a}^{\,Q_C}.
\]
The inner product $\langle\cdot,\cdot\rangle_{A_a}$ is the continuous
extension of the inner product induced by $Q_C$ on the quotient.
This definition is meaningful for every $a>0$.
To specify the value on $\Gamma$, we use
\[
 \ell^2(\Gamma,m):=
 \left\{(c_\gamma)_{\gamma\in\Gamma}:
              \sum_\gamma m_\gamma|c_\gamma|^2<\infty\right\},
 \qquad m=(m_\gamma)_{\gamma\in\Gamma},
\]
equipped with the weighted $\ell^2$ inner product determined by $m$. 
By \eqref{eq_104}, the map
\[
 u\longmapsto(\widehat u(\gamma))_{\gamma\in\Gamma},
 \qquad u\in C_c^\infty(-a,a),
\]
passes to the quotient by the null space of $Q_C$ and extends by
completion to an isometric embedding
\[
 \Samp_a:\Ha\longrightarrow\ell^2(\Gamma,m).
\]
We call this the sampling map.

\subsection{Decomposition by periods and Vandermonde matrices}

Fix $a>0$ and extend $u\in L^2(-a,a)$ by zero.
For a fixed $t\in(0,h)$, we call the sequence of points
$\{t+kh:k\in\Z\}$ a fiber.
For almost every $t\in(0,h)$, we define
\[
 L_j u(t):=e^{i\beta_jt}
             \sum_{k\in\Z}\zeta_j^k u(t+kh).
\]
The sum is finite.
The change of variables $x=t+kh$ and the identity
$\zeta_j=e^{i\beta_jh}$ give
\[
 \widehat u(\beta_j+2\pi n/h)
 =\int_0^h L_j u(t)e^{2\pi int/h}\,dt.
\]
Parseval's identity for Fourier series yields
\begin{equation} \label{eq_401}
 \begin{split}
 Q_C(u)&=h\int_0^h\sum_{j=1}^d m_j
             \left|\sum_{k\in\Z}\zeta_j^ku(t+kh)\right|^2\,dt,\\
 \|u\|_{L^2}^2&=\int_0^h
             \sum_{k\in\Z}|u(t+kh)|^2\,dt.
 \end{split}
\end{equation}
Thus the sum over zeros in \eqref{eq_104} converges for compactly supported
$L^2$ functions as well.
This expression agrees with $Q_C$ on test functions, 
and both forms are bounded on each fixed interval $(-a,a)$.
It therefore agrees with the extension to $L^2(-a,a)$ defined in
Section~\ref{section_3}.
For each $t\in(0,h)$, put
\[
 I_a(t):=\{k\in\Z:-a<t+kh<a\},\qquad
 r(t):=\#I_a(t).
\]
When $r(t)\ge1$, set $k_0(t):=\min I_a(t)$.
Then $I_a(t)$ consists of the consecutive integers
$k_0(t),\ldots,k_0(t)+r(t)-1$.
Define the column vector
\[
 \mathbf u(t):=
 \bigl(u(t+(k_0(t)+\ell)h)\bigr)_{0\le\ell<r(t)}
 \in\C^{r(t)}.
\]
For $r(t)=0$, we take $\C^0=\{0\}$.
The second identity in \eqref{eq_401} shows that
$u\mapsto\mathbf u$ gives an isometric isomorphism
\[
\begin{aligned}
 L^2(-a,a)
 &\simeq \int_{(0,h)}^{\oplus}\C^{r(t)}\,dt\\
 &:=\left\{
 \mathbf u(t)\in\C^{r(t)}
 \ \text{a.e. }t\in(0,h):
 \int_0^h\|\mathbf u(t)\|_{\C^{r(t)}}^2\,dt<\infty
 \right\},
\end{aligned}
\]
where the direct integral is taken with respect to Lebesgue measure on $(0,h)$. 
Conversely, assigning each component to the corresponding points
$t+(k_0(t)+\ell)h$ recovers a function $u\in L^2(-a,a)$.
We call this the fiber decomposition of $L^2(-a,a)$.
For $r\ge1$ and $k_0\in\Z$, define the Vandermonde matrix
\begin{equation} \label{eq_402}
 V_r(k_0):=
 (\sqrt{m_j}\,\zeta_j^{k_0+\ell})_{
 \substack{1\le j\le d\\0\le\ell<r}}.
\end{equation}
Since
\[
 V_r(k_0)=
 \operatorname{diag}(\zeta_1^{k_0},\ldots,\zeta_d^{k_0})V_r(0),
\]
the matrix $V_r(k_0)^\ast V_r(k_0)$ is independent of $k_0$.
Moreover, the columns of $V_r(k_0)$ are linearly independent for $r\le d$.
In this notation, \eqref{eq_401} becomes
\[
 Q_C(u)=h\int_0^h
   \|V_{r(t)}(k_0(t))\mathbf u(t)\|_{\C^d}^2\,dt,
\]
where the integrand is taken to be zero when $r(t)=0$.

\begin{theorem} \label{thm_401}
There exist constants $c_\ast,C_\ast>0$, independent of $a$, such that,
for $0<a\le a_\ast$,
\begin{equation} \label{eq_403}
 c_\ast\|u\|_{L^2}^2
 \le Q_C(u)\le C_\ast\|u\|_{L^2}^2
 \qquad(u\in L^2(-a,a)).
\end{equation}
Consequently, $N_a=0$ and $\Ha=L^2(-a,a)$ as sets.
At the endpoint $a=a_\ast$, the map
\begin{equation} \label{eq_404}
 \Samp_\ast:=\Samp_{a_\ast}:
 \mathcal H(A_{a_\ast})
 \overset{\sim}{\longrightarrow}\ell^2(\Gamma,m)
\end{equation}
is an isometric isomorphism.
\end{theorem}

\begin{proof}
The upper bound in \eqref{eq_403} follows uniformly for
$a\le a_\ast$ from \eqref{eq_307}.
To prove the lower bound, note that $2a\le dh$ implies $r(t)\le d$
for almost every $t\in(0,h)$.
For $r\le d$, the columns of $V_r(0)$ are linearly independent, so
$V_r(0)^\ast V_r(0)$ is positive definite.
Thus, with
\[
 c_\ast:=h\min_{1\le r\le d}
 \lambda_{\min}\bigl(V_r(0)^\ast V_r(0)\bigr)>0,
\]
\eqref{eq_401} gives $Q_C(u)\ge c_\ast\|u\|_{L^2}^2$, 
where $\lambda_{\min}(H)$ denotes the smallest eigenvalue of a Hermitian matrix $H$. 
The density of $C_c^\infty(-a,a)$ in $L^2(-a,a)$ then gives
$\Ha=L^2(-a,a)$ as sets.

We next prove surjectivity of the sampling map $\Samp_\ast$.
Since $2a_\ast=dh$, we have $r(t)=d$ for almost every $t\in(0,h)$,
and $V_d(k_0(t))$ is invertible.
Take any $c=(c_\gamma)\in\ell^2(\Gamma,m)$.
For each $j$, the Riesz--Fischer theorem gives $f_j\in L^2(0,h)$ such that
\[
 c_{\beta_j+2\pi n/h}
 =\int_0^h f_j(t)e^{2\pi int/h}\,dt
 \qquad(n\in\Z).
\]
For almost every $t\in(0,h)$, define $\mathbf u(t)$ by solving
\[
 V_d(k_0(t))\mathbf u(t)=
 \begin{pmatrix}
 \sqrt{m_1}e^{-i\beta_1t}f_1(t)\\
 \vdots\\
 \sqrt{m_d}e^{-i\beta_dt}f_d(t)
 \end{pmatrix}.
\]
The norms of $V_d(k_0(t))^{-1}$ are uniformly bounded.
The fiber decomposition therefore gives $u\in L^2(-a_\ast,a_\ast)$.
Since $L_j u=f_j$, we have $\Samp_\ast u=c$.
The isometry follows from the definition of $\Samp_\ast$.
\end{proof}

\begin{proposition} \label{prop_402}
If $a>a_\ast$, then $N_a\ne0$.
\end{proposition}
\begin{proof}
Choose $0<\eta<a-a_\ast$ and a nonzero function
$\psi\in C_c^\infty(-\eta,\eta)$.
By \eqref{eq_204}, $\Psf$ is a finite exponential sum of
exponential type $a_\ast$.
Thus $\Psf(z)\widehat\psi(z)$ is the Fourier transform of some
$u\in C_c^\infty(-a,a)$.
This function $u$ is nonzero, and $\widehat u$ vanishes on $\Gamma$.
Hence $u\in N_a\setminus\{0\}$.
\end{proof}

\subsection{The spectrum and lowest eigenvalue of $A_a$}

We now describe the spectrum of the bounded operator $A_a$ on
$L^2(-a,a)$.
For $r\ge1$, define the Toeplitz matrix
\begin{equation} \label{eq_405}
 H_r:=(p_{\ell-k})_{0\le k,\ell<r}
     =V_r(0)^\ast V_r(0).
\end{equation}
The equality follows from $p_n=\sum_jm_j\zeta_j^n$.
Since $p_n$ is real and even in $n$, $H_r$ is real symmetric.

\begin{proposition} \label{prop_403}
Write $L=2a/h$.
If $L=n\in\N$, then
\[
 \sigma(A_a)=h\sigma(H_n).
\]
If $L=n+\eta$, where $n\ge0$ and $0<\eta<1$, then
\[
 \sigma(A_a)=h\sigma(H_n)\cup h\sigma(H_{n+1}),
\]
with the $H_n$ term omitted when $n=0$.
All these eigenvalues have infinite multiplicity, and $A_a$ has no
other spectrum.
In particular,
\begin{equation} \label{eq_406}
 \lambda_a:=\inf\sigma(A_a)
       =h\lambda_{\min}(H_{\lceil2a/h\rceil}).
\end{equation}
Moreover, $\lambda_a$ is a nonincreasing, left-continuous step function.
It is positive for $0<a\le a_\ast$ and zero for $a>a_\ast$.
\end{proposition}
\begin{proof}
Polarizing \eqref{eq_401} and using \eqref{eq_405}, we obtain
\[
 Q_C(u,v)=h\int_0^h
 \langle H_{r(t)}\mathbf u(t),
         \mathbf v(t)\rangle_{\C^{r(t)}}\,dt,
\]
where the integrand is zero when $r(t)=0$.
Thus, under the fiber decomposition of $L^2(-a,a)$, the operator $A_a$
is represented by multiplication by the matrix $hH_{r(t)}$.
The part where $r(t)=0$ is zero-dimensional and may be omitted.

If $L=n\in\N$, then $r(t)=n$ for almost every $t\in(0,h)$.
If $L=n+\eta$ with $n\ge0$ and $0<\eta<1$, then $r(t)$ is either
$n$ or $n+1$.
Since
\[
 \int_0^h r(t)\,dt=2a=(n+\eta)h,
\]
the set where $r(t)=n+1$ has measure $\eta h$, and the set where
$r(t)=n$ has measure $(1-\eta)h$.
Consequently, $A_a$ is unitarily equivalent to the direct sum, over those
$r\ge1$ for which $\{t:r(t)=r\}$ has positive measure, of multiplication
by the constant matrix $hH_r$ on the corresponding $L^2$ space.
The zero-dimensional part where $r(t)=0$ contributes no spectrum.
This proves the stated formula for $\sigma(A_a)$.
Each matrix eigenvector can be multiplied by an arbitrary $L^2$ function
on the corresponding set, so every eigenvalue has infinite multiplicity.
Since $H_r$ is a principal submatrix of $H_{r+1}$,
$\lambda_{\min}(H_{r+1})\le\lambda_{\min}(H_r)$.
This proves \eqref{eq_406}, together with monotonicity and
left continuity of $\lambda_a$.
For $r\le d$, we have $H_r>0$; for $r>d$,
$\operatorname{rank}H_r=d<r$, so the smallest eigenvalue is zero.
Thus $\lambda_a>0$ for $a\le a_\ast$ and $\lambda_a=0$ for $a>a_\ast$.
\end{proof}

\begin{corollary} \label{cor_404}
For every $a>0$, the value $\lambda_a$ is an eigenvalue and is attained
by a nonzero test function $u\in C_c^\infty(-a,a)$.
In particular, for $0<2a\le h$,
\begin{equation} \label{eq_407}
 A_a=hp_0I,\qquad \lambda_a=hp_0.
\end{equation}
In this range, every nonzero even or odd function is an eigenfunction
corresponding to the lowest eigenvalue.
\end{corollary}
\begin{proof}
By Proposition~\ref{prop_403}, $\lambda_a$ is an eigenvalue
of a matrix $hH_r$ corresponding to a set of $t$ of positive measure.
Choose an open interval $J\subset(0,h)$ on which $r(t)=r$ and
$k_0(t)=k_0$ are constant, and let $v=(v_0,\ldots,v_{r-1})$ be a
nonzero eigenvector of $H_r$ for $\lambda_a/h$.
For a nonzero $\psi\in C_c^\infty(J)$, set
$u(t+(k_0+\ell)h)=v_\ell\psi(t)$ for $t\in J$ and $0\le\ell<r$,
and set $u=0$ elsewhere.
Then $u\in C_c^\infty(-a,a)$ and $A_au=\lambda_a u$.
Equation~\eqref{eq_407} follows from \eqref{eq_306}.
\end{proof}

\begin{example} \label{ex_405}
Consider the nonsingular projective model of the curve
$y^2+y=x^3+x+1$ over $\mathbb F_2$.
It has genus one, no affine $\mathbb F_2$-points, and one point at
infinity, so $N_1=1$.
Thus $P_C(T)=1-2T+2T^2$, $\omega_\nu=1\pm i$, 
and $d=2$, $a_\ast=h=\log2$, $p_0=2$, $p_1=\sqrt2$.
In this case, we have 
\[
 \lambda_a=
 \begin{cases}
 2h,&0<a\le h/2,\\
 (2-\sqrt2)h,&h/2<a\le h,\\
 0,&a>h.
 \end{cases}
\]
In particular, the continuity in \cite[Theorem~1.3]{Su26a} does not
extend directly to function fields.
More generally, $\lambda_{a_\ast}>0$, whereas $\lambda_a=0$ for
$a>a_\ast$, so right continuity fails at $a=a_\ast$.
\end{example}

\section{Local differentiation operators} \label{section_5}

We investigate the function field counterpart of the question addressed
in \cite[Theorem~1.5]{Su26a}.
Unless stated otherwise, we assume $0<a\le a_\ast$.
In this range, $A_a$ is bounded and invertible.

\subsection{de Branges spaces and Riesz vectors}
\label{section_4-kernel}

We show that the space of Fourier transforms of functions of $\Ha$
is a de Branges space of entire functions.
For the basic theory of de Branges spaces, see \cite{dB68};
a brief account is also given in \cite[Section~3.1]{Su23}.

Let $\Hhat:=\mathcal F\Ha$, with the inner product defined by
\[
 \langle\widehat u,\widehat f\rangle_{\widehat A_a}
 :=\langle u,f\rangle_{A_a}.
\]
As a set, $\Hhat=\PW_a$, and its norm is equivalent to the usual $L^2$ norm.
Point evaluation is continuous by \eqref{eq_403} and the
Cauchy--Schwarz inequality.
The conjugation $F\mapsto F^\sharp$ preserves $\PW_a$.
For $w\in\C\setminus\R$ and $F\in\PW_a$ with $F(w)=0$,
the transformation
$F(z)\mapsto((z-\bar w)/(z-w))F(z)$
also preserves $\PW_a$.
Both transformations preserve absolute values at the real sampling
points $\Gamma$.
By \eqref{eq_104}, both transformations preserve the norm of $\Hhat$.
This verifies the three axioms of a de Branges space.
There are no common zeros, since at each point one can choose a test
function whose Fourier transform does not vanish there.

Following the convention of \cite[Section~6.3]{Su26a}, define
\begin{equation} \label{eq_501}
 \langle u,u_{z,a}\rangle_{A_a}=\widehat u(\bar z)
       \qquad(u\in\Ha).
\end{equation}
With $e_z(x):=e^{-izx}$, we have $u_{z,a}=A_a^{-1}e_z$, and
$z\mapsto u_{z,a}$ is an entire $\Ha$-valued function.
The fiber decomposition of Section~\ref{section_4} allows us to compute
this vector using finite matrices.
When $r(t)\ge1$,
\begin{equation} \label{eq_502}
 u_{z,a}(t+(k_0(t)+\ell)h)
 =\frac1h\sum_{k=0}^{r(t)-1}
  \bigl(H_{r(t)}^{-1}\bigr)_{\ell k}\,
  e^{-iz(t+(k_0(t)+k)h)},
 \qquad 0\le\ell<r(t).
\end{equation}
Thus, for fixed $a\le a_\ast$, the vector $u_{z,a}$ can be computed
explicitly from finitely many $p_n$.
The reproducing kernel of $\Hhat$ is
\begin{equation} \label{eq_503}
 J_a(w,z):=\langle u_{\bar w,a},u_{\bar z,a}\rangle_{A_a}
       =\widehat{u_{\bar w,a}}(z),\qquad
 J_a(w,w)=\|u_{\bar w,a}\|_{A_a}^2>0,
\end{equation}
as can be checked directly.
In particular, $J_a(-i,\cdot)=\widehat{u_{i,a}}$ and
$J_a(i,\cdot)=\widehat{u_{-i,a}}$.
On the other hand, the general theory of de Branges spaces provides
a generating function $E_a$ for $\Hhat$ with no real zeros, for which
\begin{equation} \label{eq_504}
 J_a(w,z)=\frac{\overline{E_a(w)}E_a(z)-
       \overline{E_a^\sharp(w)}E_a^\sharp(z)}
       {2\pi i(\bar w-z)}.
\end{equation}

\subsection{The local differentiation operator}

We consider the densely defined differentiation operator $\D_a$ on $\Ha$ given by
\begin{equation} \label{eq_505}
 \Dom(\D_a):=C_c^\infty(-a,a)\subset\Ha,\qquad
 \D_a u:=iu'.
\end{equation}
Its symmetry follows from \eqref{eq_104} and the fact that the sampling
points in $\Gamma$ are real.

\begin{proposition} \label{prop_501}
The closure and adjoint of the symmetric operator $\D_a$ are given by
\begin{align}
 \Dom(\overline\D_a)&=H_0^1(-a,a),&
 \overline\D_a u&=iu',
       \label{eq_506}\\
 \Dom(\D_a^\ast)&=\{f\in L^2(-a,a):A_a f\in H^1(-a,a)\},&
 \D_a^\ast f&=A_a^{-1}i(A_a f)'
       \label{eq_507}.
\end{align}
In particular,
\begin{equation} \label{eq_508}
 \ker(\D_a^\ast-z)=\C u_{z,a},\qquad
 \D_a^\ast u_{z,a}=zu_{z,a}
 \qquad(z\in\C),
\end{equation}
and the deficiency indices are $(1,1)$.
\end{proposition}
\begin{proof}
By \eqref{eq_403}, the graph norm on test functions is
equivalent to the $H^1$-norm, so the domain of the closure is
$H_0^1(-a,a)$.
Since $A_a$ is self-adjoint, the conditions
$f\in\Dom(\D_a^\ast)$ and $\D_a^\ast f=g$ are equivalent to
\[
 \langle iu',A_a f\rangle_{L^2}
 =\langle u,A_a g\rangle_{L^2}
 \qquad(u\in C_c^\infty(-a,a)).
\]
This means that $i(A_a f)'=A_a g$ in the sense of distributions.
Since $A_a$ is boundedly invertible, \eqref{eq_507} follows.
The eigenvalue equation $\D_a^\ast f=zf$ is equivalent to
$A_a f$ being a constant multiple of $e^{-izx}$.
Since $A_a$ is boundedly invertible, the corresponding
eigenspace is spanned by the nonzero vector
$u_{z,a}=A_a^{-1}e_z$.
This proves \eqref{eq_508}.
\end{proof}

Reflection, defined by $\Jref u(x):=u(-x)$, is a unitary involution on
$\Ha$.
Since $A_a$ commutes with reflection, $\Jref u_{z,a}=u_{-z,a}$ and
$\Jref\D_a\Jref=-\D_a$.
Thus we may use
\[
 u_+(a):=u_{i,a},\qquad u_-(a):=u_{-i,a},\qquad
 \rho_a^2:=J_a(i,i)=J_a(-i,-i)=\|u_{\pm i,a}\|_{A_a}^2
\]
as deficiency vectors of equal norm in describing the self-adjoint
extensions.
By von Neumann's theory \cite[Theorem~X.2 and its Corollary]{RSII},
all self-adjoint extensions of $\D_a$ are given by
\begin{equation} \label{eq_509}
 \begin{gathered}
 w_{a,\theta}:=u_{i,a}+e^{i\theta}u_{-i,a},
 \qquad 0\le\theta<2\pi,\\
 \Dom(\overline\D_{a,\theta})
 =H_0^1(-a,a)\dotplus\C w_{a,\theta},\qquad
 \overline\D_{a,\theta}
 :=\D_a^\ast|_{\Dom(\overline\D_{a,\theta})}.
 \end{gathered}
\end{equation}
Define the boundary form by
\[
 \mathcal W_a(f,g):=\langle\D_a^\ast f,g\rangle_{A_a}
                 -\langle f,\D_a^\ast g\rangle_{A_a}.
\]
For $f\in\Dom(\D_a^\ast)$, membership in
$\Dom(\overline\D_{a,\theta})$ is equivalent to
$\mathcal W_a(f,w_{a,\theta})=0$.

\subsection{The spectral function of the local differentiation operator}

Using the Riesz vectors $u_{\pm i,a}$, define
\begin{equation} \label{eq_510}
 \begin{split}
 W(a,\theta;z)&:=(z-i)\widehat{u_{i,a}}(z)
           +e^{i\theta}(z+i)\widehat{u_{-i,a}}(z)\\
 &=(z-i)J_a(-i,z)+e^{i\theta}(z+i)J_a(i,z)
=\overline{\mathcal W_a(u_{\bar z,a},w_{a,\theta})}.
 \end{split}
\end{equation}

\begin{theorem} \label{thm_502}
The self-adjoint extension $\overline\D_{a,\theta}$ has compact
resolvent.
The function $W(a,\theta;z)$ is a nonzero entire function of exponential
type at most $a$.
Its zeros are all real and simple and coincide with the spectrum of
$\overline\D_{a,\theta}$.
The domain of $\overline\D_{a,\theta}$ is invariant under reflection
if and only if $\theta=0$ or $\pi$.
The function $W(a,0;z)$ is odd and $W(a,\pi;z)$ is even.
Moreover, $\ker\overline\D_{a,0}=\C u_{0,a}$ and
$0\in\rho(\overline\D_{a,\pi})$.
\end{theorem}

\begin{proof}
If a family of vectors $f\in\Dom(\D_a^\ast)$ is bounded in the graph
norm, then both $A_a f$ and
$(A_a f)'=-iA_a\D_a^\ast f$ are bounded in $L^2(-a,a)$.
The compactness of $H^1(-a,a)\hookrightarrow L^2(-a,a)$ and
boundedness of $A_a^{-1}$ therefore imply that the inclusion of
$\Dom(\D_a^\ast)$, equipped with its graph norm, into $\Ha$ is compact.
Hence every self-adjoint extension has compact resolvent.

Since $u_{\pm i,a}\in L^2(-a,a)$, the function in \eqref{eq_510}
is entire and satisfies
$|W(a,\theta;z)|\le C_a(1+|z|)e^{a|\Im z|}$.
Also, $W(a,\theta;i)=2ie^{i\theta}\rho_a^2\ne0$.
By \eqref{eq_510} and the boundary condition, if $W(a,\theta;z)=0$,
then the nonzero vector $u_{\bar z,a}$ is an eigenvector of the
extension for $\bar z$.
Self-adjointness implies that $z$ is real.
Conversely, by \eqref{eq_508}, an eigenvector for a real eigenvalue
$z$ is a multiple of $u_{z,a}$, so $W(a,\theta;z)=0$.
Compactness of the resolvent shows that there is no other spectrum.

To prove simplicity of the zeros as zeros of an entire function, we use
the reproducing kernel.
Put $p_a(z):=(z+i)J_a(i,z)$.
Then $p_a^\sharp(z)=(z-i)J_a(-i,z)$, and \eqref{eq_504} gives
\[
 p_a(z)=\frac{i}{2\pi}
   \bigl(\overline{E_a(i)}E_a(z)
   -\overline{E_a^\sharp(i)}E_a^\sharp(z)\bigr),
 \quad
 |E_a(i)|^2-|E_a^\sharp(i)|^2=4\pi\rho_a^2.
\]
Substitution into the kernel formula yields
\begin{equation} \label{eq_511}
 \frac{\overline{p_a(w)}p_a(z)-
       \overline{p_a^\sharp(w)}p_a^\sharp(z)}
      {2\pi i(\bar w-z)}
       =\frac{\rho_a^2}{\pi}J_a(w,z).
\end{equation}
Define real entire functions $s_{a,\theta}$ and $c_{a,\theta}$ by
$e^{i\theta/2}p_a(z)=s_{a,\theta}(z)+ic_{a,\theta}(z)$.
Multiplication by this phase does not change \eqref{eq_511},
and $W(a,\theta;z)=2e^{i\theta/2}s_{a,\theta}(z)$.
Taking the diagonal limit at a real zero $t$, we obtain
\[
 0<\rho_a^2J_a(t,t)=c_{a,\theta}(t)s_{a,\theta}'(t).
\]
Thus $s_{a,\theta}'(t)\ne0$, proving simplicity.

Reflection sends $w_{a,\theta}$ to $e^{i\theta}w_{a,-\theta}$.
Since distinct phases give distinct domains, invariance under reflection
is equivalent to $\theta=0$ or $\pi$.
The identity $J_a(-i,z)=J_a(i,-z)$ shows that $W(a,0)$ is odd and
$W(a,\pi)$ is even.
Since $u_{0,a}$ and $w_{a,0}$ are even, whereas $\D_a^\ast w_{a,0}$
is odd, we have $\mathcal W_a(u_{0,a},w_{a,0})=0$.
Thus zero is an eigenvalue for phase $0$, and \eqref{eq_508}
shows that its eigenspace is one-dimensional.
For phase $\pi$, the function $W$ is even and has only simple zeros,
so it cannot vanish at the origin.
The description of the spectrum then gives
$0\in\rho(\overline\D_{a,\pi})$.
\end{proof}

For any $a>a_\ast$, a similar local theory can be constructed by
choosing $\lambda<\inf\sigma(A_a)$ and using
$T_{a,\lambda}=A_a-\lambda I$.
The operator $T_{a,\lambda}$ is bounded and positive definite, and
$\mathcal H(T_{a,\lambda})=L^2(-a,a)$ as sets.
Replacing $A_a$ by $T_{a,\lambda}$ in the proof of
Proposition~\ref{prop_501} gives deficiency indices $(1,1)$
and self-adjoint extensions.
The argument for de Branges spaces also applies, since
$Q_C(u)-\lambda\|u\|_2^2$ is unchanged by moving a nonreal zero in a
way that preserves absolute values on the real axis.
This corresponds to the construction on arbitrary finite intervals
$(-a,a)$ in \cite[Theorem~1.5]{Su26a}.
The construction at $a=a_\ast$ in the next section and the limits of
compressions, however, concern $\lambda=0$ and $0<a\le a_\ast$.
\medskip

In the next section, we show that the self-adjoint extension
$\overline\D_{a,\theta_\ast}$ corresponding to the distinguished
phase \eqref{eq_105} has spectrum $\Gamma$ at $a=a_\ast$.

\section{The spectrum of the local differentiation operator at the endpoint}
\label{section_6}

We investigate the function field counterparts of
\cite[Theorems~1.1 and 1.2]{Su26b} and, in particular, prove
Theorems~\ref{thm_1} and \ref{thm_2}.

Write $\Hstar:=\mathcal H(A_{a_\ast})$, and denote its inner product
and norm by $\langle\cdot,\cdot\rangle_\ast$ and $\|\cdot\|_\ast$,
respectively.
By Theorem~\ref{thm_401},
$\widehat{\Hstar}=(\PW_{a_\ast},\widehat Q_C)$ is a de Branges space,
and sampling its elements on $\Gamma$ gives an isometric isomorphism
onto $\ell^2(\Gamma,m)$, 
where $\widehat Q_C(\widehat u,\widehat f):=Q_C(u,f)$.

\subsection{An explicit formula for the reproducing kernel at the endpoint}
\label{section_6_1}
For each $\gamma\in\Gamma$, we define
\[
 d_\gamma(z):=\frac{\Psf(z)}{\Psf'(\gamma)(z-\gamma)}.
\]
This function belongs to $\PW_{a_\ast}$.
Indeed, the pole is removable, its exponential type is at most
$a_\ast$, and it is $O(1/|z|)$ at infinity on the real axis, so the
Paley--Wiener theorem applies.
Since $d_\gamma(\eta)=\delta_{\gamma\eta}$ and the sampling map is
surjective, $\{d_\gamma/\sqrt{m_\gamma}\}$ is an orthonormal basis
of $\widehat{\Hstar}$.
Consequently, its reproducing kernel is
\begin{equation} \label{eq_601}
 J_\ast(w,z)=\sum_{\gamma\in\Gamma}
         \frac{\overline{d_\gamma(w)}d_\gamma(z)}{m_\gamma},
 \qquad
 J_\ast(\gamma,\gamma)=\frac1{m_\gamma}.
\end{equation}
The series converges uniformly on compact sets in each variable.
By \eqref{eq_204},
$\Psf'(z+2\pi/h)=(-1)^d\Psf'(z)$, so $\Psf'(\gamma)^2$ takes only
finitely many nonzero values.
For every compact set $K\subset\C$, therefore,
$|d_\gamma(z)|\ll_K(1+|\gamma|)^{-1}$ on $K$, which gives absolute
uniform convergence on compact sets.

To compute the kernel, we introduce the auxiliary entire functions
\begin{equation} \label{eq_602}
\aligned
 B_\ast(z)&:=\Psf(z)\sum_{j=1}^d
          \frac{\pi h}{2m_j\Psf'(\beta_j)^2}
                  \cot\frac{h(z-\beta_j)}2,\\
 E_\ast(z)&:=\Psf(z)+iB_\ast(z).
\endaligned
\end{equation}
The poles of each cotangent are canceled by $\Psf$.
Thus $B_\ast$ is real entire, and
\begin{equation} \label{eq_603}
 B_\ast(\gamma)=\frac{\pi}{m_\gamma\Psf'(\gamma)}.
\end{equation}
For $\Im z>0$, we have $\Im\bigl(\cot(h(z-\beta_j)/2)\bigr)<0$, so
\[
 |E_\ast(z)|^2-|E_\ast^\sharp(z)|^2
 =4|\Psf(z)|^2
 \Im\left(-\frac{B_\ast(z)}{\Psf(z)}\right)>0.
\]
Moreover, \eqref{eq_603} shows that $\Psf$ and $B_\ast$
have no common real zeros.
Hence $E_\ast$ is a Hermite--Biehler function without real zeros.

\begin{lemma} \label{lem_601}
The de Branges space associated with $E_\ast$ coincides isometrically
with $\widehat{\Hstar}$.
More precisely,
\begin{equation} \label{eq_604}
 \mathcal H(E_\ast)=(\PW_{a_\ast},\widehat Q_C),
 \qquad
 J_\ast(w,z)=\frac{\overline{\Psf(w)}B_\ast(z)-
                   \overline{B_\ast(w)}\Psf(z)}
                  {\pi(\bar w-z)}.
\end{equation}
Furthermore, $B_\ast(-z)=-\varepsilon B_\ast(z)$ and $B_\ast(i)\ne0$.
\end{lemma}
\begin{proof}
We have $\Psf'(\beta_j+2\pi n/h)^2=\Psf'(\beta_j)^2$.
Write $r_\ast(z):=-B_\ast(z)/\Psf(z)$.
Taking the difference of the partial fraction expansions of the cotangent
\cite[(4.22.3)]{NIST}, we obtain
\[
 \frac{r_\ast(z)-r_\ast(\bar w)}{\pi(z-\bar w)}
 =\sum_{\gamma\in\Gamma}
   \frac{1}{m_\gamma\Psf'(\gamma)^2(z-\gamma)(\bar w-\gamma)}.
\]
The series of differences converges absolutely and uniformly on compact
sets avoiding the poles.
Multiplication by $\overline{\Psf(w)}\Psf(z)$ shows that
\eqref{eq_601} and \eqref{eq_604} agree.
The right-hand side of \eqref{eq_604} is also the reproducing kernel
of $\mathcal H(E_\ast)$, so equality of the kernels gives the
isometric identification of the Hilbert spaces.
The quantities $m_j\Psf'(\beta_j)^2$ agree for residue classes paired
by reflection.
Thus $r_\ast$ is odd, and $B_\ast(-z)=-\varepsilon B_\ast(z)$.
Since $\Psf(i)\ne0$ and $\Im r_\ast(i)>0$, we also have
$B_\ast(i)\ne0$.
\end{proof}

\subsection{Convergence through orthogonal projections}
\label{sec:W-limit}

By extension by zero, the spaces $\Ha$ form an increasing family of
closed subspaces of $\Hstar$.
Let $P_a^Q:\Hstar\to\Hstar$ be the $Q_C$-orthogonal projection onto
this subspace.
The density of test functions gives
\[
 \overline{\bigcup_{0<a<a_\ast}\Ha}=\Hstar,
\]
so $P_a^Q\to I$ strongly.
Writing $u_{z,\ast}:=u_{z,a_\ast}$ 
for the Riesz vector at the endpoint, uniqueness in \eqref{eq_501}
gives $u_{z,a}=P_a^Q u_{z,\ast}$. 
Consequently, we have 
\begin{equation} \label{eq_605}
 u_{z,a}~\to~u_{z,\ast}\quad\text{in }Q_C,\qquad
 J_a(w,z)~\to~J_\ast(w,z)\quad(a\uparrow a_\ast).
\end{equation}
The Riesz vector $u_{z,\ast}=A_{a_\ast}^{-1}e_z$ depends continuously
on $z$ in norm, and the image of a compact set is compact in the norm
topology.
Strong convergence of uniformly bounded operators is uniform on such
sets.
Thus \eqref{eq_605} holds uniformly on compact sets of $z$ and of
$(w,z)$, respectively.

\begin{theorem} \label{thm_602}
For the phase in \eqref{eq_105}, we have 
\begin{equation} \label{eq_606}
 \lim_{a\uparrow a_\ast}W(a,\theta_\ast;z)
 =\frac{2B_\ast(i)}{\pi}\Psf(z), 
\end{equation}
where the convergence is uniform on every compact subset of $\C$.
At the endpoint, the identity
$W(a_\ast,\theta_\ast;z)=2B_\ast(i)\Psf(z)/\pi$ holds, and the
spectrum of the self-adjoint extension $\overline\D_{a_\ast,\theta_\ast}$ coincides with
$\Gamma$.
\end{theorem}

\begin{proof}
Substituting $w=-i,i$ into \eqref{eq_604} and using reality and
parity, we obtain
\begin{align*}
 (z-i)J_\ast(-i,z)
 &=\frac{B_\ast(i)\Psf(z)-\Psf(i)B_\ast(z)}{\pi},\\
 (z+i)J_\ast(i,z)
 &=\frac{-\varepsilon B_\ast(i)\Psf(z)
           -\varepsilon\Psf(i)B_\ast(z)}{\pi}.
\end{align*}
We have $e^{i\theta_\ast}=-\varepsilon$, so the terms involving $B_\ast(z)$ cancel.
This proves the identity at the endpoint. 
The locally uniform convergence in \eqref{eq_606} follows from
\eqref{eq_510} and \eqref{eq_605}.
The zeros of $\Psf$ are precisely $\Gamma$, and
Theorem~\ref{thm_502} applies at the endpoint, proving the
assertion about the spectrum.
\end{proof}

\subsection{Proof of Theorem~\ref{thm_2}}

By Theorem~\ref{thm_502}, we have $W(a,\theta_\ast;i)\ne0$.
The limit in Theorem~\ref{thm_602} at $z=i$ is $2B_\ast(i)\Psf(i)/\pi\ne0$.
Therefore, taking the ratio yields \eqref{eq_107}. 
For a curve, if all zeros are simple, then there is no central zero
because $A_C$ is even, and $\Psf$ differs from $A_C$ by a nonzero
constant factor.
Thus $\theta_\ast=\pi$, and the additional assertion following
\eqref{eq_107} also holds.
\hfill $\Box$

\subsection{The differentiation operator at the endpoint}
\label{section_6-compression}

We prepare the proof of monotonicity in Theorem~\ref{thm_1}.
For simplicity, write $D_a:=\overline\D_{a,\theta_\ast}$ and
$D_\ast:=D_{a_\ast}$.
Let $\iota_a:\Ha\hookrightarrow\Hstar$ be the isometric embedding
given by extension by zero, and denote its adjoint by $\iota_a^\ast$.
By \eqref{eq_403} and the argument preceding \eqref{eq_605}, we have 
\begin{equation} \label{eq_607}
 \iota_a^\ast\iota_a=I_{\Ha},\qquad
 P_a^Q=\iota_a\iota_a^\ast,\qquad
 \iota_a\Ha\subset\iota_b\mathcal H(A_b)\ (a<b\le a_\ast),
\end{equation}
and $P_a^Q\to I_{\Hstar}$ strongly.

\begin{lemma} \label{lem_603}
The vectors $e_{\gamma,\ast}:=\sqrt{m_\gamma}\,u_{\gamma,\ast}$
$(\gamma\in\Gamma)$ form an orthonormal basis of $\Hstar$, and
$D_\ast e_{\gamma,\ast}=\gamma e_{\gamma,\ast}$.
Moreover,
\begin{equation} \label{eq_608}
 \begin{split}
 \Dom(D_\ast)
 &=\bigl\{f\in\Hstar:
   \sum_{\gamma\in\Gamma}m_\gamma\gamma^2
             |\widehat f(\gamma)|^2<\infty\bigr\},\\
 \widehat{D_\ast f}(\gamma)
 &=\gamma\widehat f(\gamma)
 \qquad(f\in\Dom(D_\ast)).
 \end{split}
\end{equation}
\end{lemma}
\begin{proof}
By \eqref{eq_501} and \eqref{eq_601}, we have
$\widehat{e_{\gamma,\ast}}=d_\gamma/\sqrt{m_\gamma}$.
Since $\{d_\gamma/\sqrt{m_\gamma}\}_{\gamma\in\Gamma}$
is an orthonormal basis of $\widehat{\Hstar}$ by
Section~\ref{section_6_1},
the vectors $\{e_{\gamma,\ast}\}_{\gamma\in\Gamma}$
form an orthonormal basis of $\Hstar$.
By \eqref{eq_508}, each $e_{\gamma,\ast}$
is an eigenvector of $\D_{a_\ast}^\ast$
with eigenvalue $\gamma$.
Theorem~\ref{thm_602} gives
$W(a_\ast,\theta_\ast;\gamma)=0$.
Hence Theorem~\ref{thm_502} shows that
$e_{\gamma,\ast}\in\Dom(D_\ast)$ and
$D_\ast e_{\gamma,\ast}=\gamma e_{\gamma,\ast}$.
For $f\in\Hstar$, write
$f=\sum_{\gamma\in\Gamma}c_\gamma e_{\gamma,\ast}$,
where $c_\gamma=\sqrt{m_\gamma}\widehat f(\gamma)$.
Since $D_\ast$ is self-adjoint, the spectral theorem
shows that $f$ belongs to $\Dom(D_\ast)$ if and only if
$\sum_{\gamma\in\Gamma}\gamma^2|c_\gamma|^2<\infty$.
For such $f$, we also have
$D_\ast f=\sum_{\gamma\in\Gamma}
\gamma c_\gamma e_{\gamma,\ast}$.
This proves \eqref{eq_608}.
\end{proof}

For $u\in C_c^\infty(-a,a)$, we have
\begin{equation} \label{eq_609}
 \iota_a u\in\Dom(D_\ast),\qquad
 D_\ast\iota_a u=\iota_a\D_a u.
\end{equation}
Set
\[
 V_a:=\Ha\ominus\ker D_a,\qquad
 V_\ast:=\Hstar\ominus\ker D_\ast.
\]
If $0\in\Gamma$, then \eqref{eq_508} and the Riesz representation
give $V_a=\{f\in\Ha:\widehat f(0)=0\}$ and
$V_\ast=\{f\in\Hstar:\widehat f(0)=0\}$.
Thus, whether or not there is a central zero,
\[
 \kappa_a:=\iota_a|_{V_a}:V_a\to V_\ast
\]
is an isometric embedding.
The subspaces $V_a$ and $V_\ast$ reduce $D_a$ and $D_\ast$,
respectively.
Define
\[
 \mathsf C_\ast:=(D_\ast|_{V_\ast})^{-1},\qquad
 \mathsf C_a:=\kappa_a^\ast\mathsf C_\ast\kappa_a.
\]
These are bounded compact self-adjoint operators.
We now show that $\mathsf C_a=(D_a|_{V_a})^{-1}$.

\begin{proposition} \label{prop_604}
Suppose $0\notin\Gamma$.
Then $\theta_\ast=\pi$, $V_a=\Ha$, $V_\ast=\Hstar$, and
$\kappa_a=\iota_a$.
For every $0<a\le a_\ast$,
\begin{equation} \label{eq_610}
 D_a^{-1}=\mathsf C_a=\kappa_a^\ast\mathsf C_\ast\kappa_a,
\end{equation}
with equality of domains.
\end{proposition}
\begin{proof}
By Lemma~\ref{lem_603}, $\mathsf C_\ast$ acts by
multiplication by $1/\gamma$ on the eigenvector basis of $D_\ast$.
Since $\Gamma$ is a finite union of arithmetic progressions and
$0\notin\Gamma$, we have
$\sum_{\gamma\in\Gamma}\gamma^{-2}<\infty$.
Thus $\mathsf C_\ast$ is a bounded Hilbert--Schmidt self-adjoint
operator, and $\mathsf C_a$ is compact and self-adjoint.
For every $u\in C_c^\infty(-a,a)$,
\eqref{eq_609} gives
\begin{equation} \label{eq_611}
 \mathsf C_a\D_a u
 =\iota_a^\ast\mathsf C_\ast D_\ast\iota_a u=u.
\end{equation}
If $f\in\ker\mathsf C_a$, then, for every test function $u$,
\[
 \langle u,f\rangle_{A_a}
 =\langle\mathsf C_a\D_a u,f\rangle_{A_a}
 =\langle\D_a u,\mathsf C_a f\rangle_{A_a}=0.
\]
Since the test functions are dense in $\Ha$, we obtain $f=0$.
Thus $\ker\mathsf C_a=\{0\}$.
As $\mathsf C_a$ is self-adjoint, its range is dense in $\Ha$. 
Define $L_a:=\mathsf C_a^{-1}$ on
$\Dom(L_a)=\ran\mathsf C_a$.
This operator is self-adjoint.
Indeed, if $h\in\Dom(L_a^\ast)$ and $k=L_a^\ast h$, then, for all
$f\in\Ha$,
\[
 \langle f,h\rangle_{A_a}
 =\langle L_a\mathsf C_a f,h\rangle_{A_a}
 =\langle\mathsf C_a f,k\rangle_{A_a}
 =\langle f,\mathsf C_a k\rangle_{A_a}.
\]
Hence $h=\mathsf C_a k\in\Dom(L_a)$ and $L_ah=k$.
Together with symmetry, this proves $L_a^\ast=L_a$.
Since $L_a$ is bijective and its inverse $\mathsf C_a$ is bounded,
$0\in\rho(L_a)$.
Equation~\eqref{eq_611} shows that $L_a$ is a
self-adjoint extension of $\D_a$.
For reflection $\Jref f(x)=f(-x)$, the identities
$\Jref D_\ast\Jref=-D_\ast$ and $\Jref\iota_a=\iota_a\Jref$ give
$\Jref\mathsf C_a\Jref=-\mathsf C_a$.
Thus $\Dom(L_a)$ is invariant under reflection.
By Theorem~\ref{thm_502}, its phase is either $0$ or $\pi$.
The extension with phase $0$ has zero as an eigenvalue, whereas
$0\in\rho(L_a)$.
Therefore the phase is $\pi$, and $L_a=D_a$.
\end{proof}

\begin{proposition} \label{prop_605}
Suppose $0\in\Gamma$.
Then $V_a=\Ha\ominus\C u_{0,a}$ and
$V_\ast=\Hstar\ominus\C u_{0,\ast}$, and
\begin{equation} \label{eq_612}
 (D_a|_{V_a})^{-1}
 =\mathsf C_a=\kappa_a^\ast\mathsf C_\ast\kappa_a
\end{equation}
with equality of domains.
\end{proposition}

\begin{proof}
Each of $D_a$ and $D_\ast$ decomposes as an orthogonal direct sum on
its zero eigenspace and the orthogonal complement.
The restricted operators have domains $\Dom(D_a)\cap V_a$ and
$\Dom(D_\ast)\cap V_\ast$, respectively.
Since zero is an isolated simple eigenvalue and both operators have
compact resolvent, these restrictions are self-adjoint and have bounded
inverses.
The operator $\mathsf C_\ast$ is Hilbert--Schmidt and self-adjoint,
acting by $1/\gamma$ on the basis of eigenvectors corresponding to
nonzero eigenvalues of $D_\ast$.
For $f\in V_a$, set
$w:=\mathsf C_a f=\kappa_a^\ast\mathsf C_\ast\kappa_a f\in V_a$.
For every $u\in C_c^\infty(-a,a)$, we have
$\widehat{\D_a u}(0)=0$.
Thus \eqref{eq_609} and self-adjointness of $D_\ast$
give
\[
 \begin{aligned}
 \langle\D_a u,w\rangle_{A_a}
 &=\langle\iota_a\D_a u,\mathsf C_\ast\kappa_a f\rangle_\ast
 =\langle D_\ast\iota_a u,\mathsf C_\ast\kappa_a f\rangle_\ast\\
 &=\langle\iota_a u,\iota_a f\rangle_\ast
 =\langle u,f\rangle_{A_a}.
 \end{aligned}
\]
Hence $w\in\Dom(\D_a^\ast)$ and $\D_a^\ast w=f$.
The vector $v:=(D_a|_{V_a})^{-1}f\in V_a$ satisfies the same equation.
By \eqref{eq_508}, $w-v\in\ker\D_a^\ast=\C u_{0,a}$.
Since $w-v$ also belongs to $V_a$, we have $w=v$.
This proves \eqref{eq_612}.
\end{proof}

\begin{remark}
Suppose $0\in\Gamma$.
We extend $\mathsf C_\ast$ by zero on $\ker D_\ast$ to obtain the
pseudoinverse $\mathsf C_\ast^0$.
However, the compression of $\mathsf C_\ast^0$ to $\Ha$
does not, in general, coincide with the local pseudoinverse. 
Indeed, let $P_{0,\ast}$ be the $Q_C$-orthogonal projection onto the
span of $u_{0,\ast}$.
For a test function,
\[
 \iota_a^\ast\mathsf C_\ast^0\iota_a\D_a u
 =u-\iota_a^\ast P_{0,\ast}\iota_a u
 =u-\frac{\widehat u(0)}{\|u_{0,\ast}\|_\ast^2}u_{0,a}.
\]
The projection onto the local zero eigenspace, however, is
$P_{0,a}u=\widehat u(0)u_{0,a}/\|u_{0,a}\|_{A_a}^2$, which
generally differs from the compressed projection.
In Proposition~\ref{prop_605}, both zero
eigenspaces are removed before taking the compression.
\end{remark}

\subsection{Proof of Theorem~\ref{thm_1}}
\label{sec:finite-minmax}

Theorem~\ref{thm_502} shows that $D_a$ has compact
resolvent and simple eigenvalues.
Theorem~\ref{thm_602} identifies the spectrum of $D_\ast$ with
$\Gamma$.
Reflection anticommutes with the extension at the distinguished phase,
so the spectrum is symmetric about the origin.
Since the local operator has compact resolvent on an infinite-dimensional
space, there are infinitely many positive and negative eigenvalues.
Zero is a simple eigenvalue precisely in the case of a central zero.

We first establish the limit.
Fix $n$ and choose $\gamma_n<R<\gamma_{n+1}$.
The zeros of $\Psf$ in $|z|<R$ are
$\pm\gamma_1,\ldots,\pm\gamma_n$, together with the origin if
$0\in\Gamma$, and they are all simple.
Choose disjoint small disks around these zeros, contained in $|z|<R$,
such that $\Psf$ does not vanish on their boundaries or on $|z|=R$.
By Theorem~\ref{thm_602} and Rouch\'e's theorem, for $a$
sufficiently close to $a_\ast$, each small disk contains exactly one
zero of $W(a,\theta_\ast;\cdot)$, and its total number of zeros in
$|z|<R$ is the same as that of $\Psf$.
There are therefore no additional zeros in $|z|<R$ outside these disks.
These zeros are real by Theorem~\ref{thm_502}.
If $0\in\Gamma$, the zero in the disk around the origin is always
$z=0$ itself.
Thus the zeros in the disks centered at positive zeros of $\Psf$ are 
$t_{a,1},\ldots,t_{a,n}$ in order.
Since the disks can be arbitrarily small, we obtain convergence of the
positive eigenvalues with their numbering preserved:
\begin{equation} \label{eq_613}
 t_{a,n}~\to~\gamma_n\qquad(a\uparrow a_\ast).
\end{equation}

For monotonicity, we use Propositions~\ref{prop_604} and
\ref{prop_605}.
In both cases,
\begin{equation} \label{eq_614}
 \mathsf C_a=\kappa_a^\ast\mathsf C_\ast\kappa_a
            =(D_a|_{V_a})^{-1}.
\end{equation}
The spaces $\kappa_a V_a$ increase with $a$, since the mean-zero
condition is preserved by extension by zero.
The positive eigenvalues of $\mathsf C_\ast$ 
are $1/\gamma_n$ in decreasing order, and those of $\mathsf C_a$ are
$\mu_{a,n}:=1/t_{a,n}$.

Let $T$ be a compact self-adjoint operator on a Hilbert space
$\mathcal H$ with at least $n$ positive eigenvalues.
Denote its $n$th positive eigenvalue by $\lambda_n^+(T)$ 
listed in decreasing order and counted with multiplicity.
The max--min principle for compact self-adjoint operators
\cite[Theorem~4.22]{Bo20} gives
\[
 \lambda_n^+(T)=
 \max_{\substack{V\subset\mathcal H\\\dim V=n}}
 \ \min_{\substack{x\in V\\\|x\|=1}}
 \langle Tx,x\rangle_{\mathcal H}.
\]
Applying \eqref{eq_614}, we obtain
\[
 \mu_{a,n}=
 \max_{\substack{V\subset\kappa_a V_a\\\dim V=n}}
 \ \min_{\substack{x\in V\\\|x\|_\ast=1}}
 \langle\mathsf C_\ast x,x\rangle_\ast.
\]
Since $\kappa_a V_a\subset\kappa_b V_b\subset V_\ast$ for
$0<a<b\le a_\ast$,
\[
 0<\mu_{a,n}\le\mu_{b,n}\le\frac1{\gamma_n},
 \qquad
 t_{a,n}\ge t_{b,n}\ge\gamma_n.
\]
Combining this with \eqref{eq_613} gives
$t_{a,n}\searrow\gamma_n$.
We have already proved $t_{a_\ast,n}=\gamma_n$ at the endpoint.
This completes the proof of Theorem~\ref{thm_1}.
\hfill $\Box$

\subsection{An alternative representation by products of eigenvalues}

Using the eigenvalues in Theorem~\ref{thm_1}, we have
\begin{equation} \label{eq_615}
 \frac{W(a,\theta_\ast;z)}{W(a,\theta_\ast;i)}
   =\left(\frac zi\right)^\chi
       \prod_{n=1}^{\infty}
       \frac{1-z^2/t_{a,n}^2}{1+1/t_{a,n}^2}, \qquad 
 \frac{\Psf(z)}{\Psf(i)}
   =\left(\frac zi\right)^\chi
       \prod_{n=1}^{\infty}
       \frac{1-z^2/\gamma_n^2}{1+1/\gamma_n^2}, 
\end{equation}
where $\chi$ is defined in \eqref{eq_102}.
These identities follow from Hadamard factorization
and the parity of the two entire functions. 
When $\chi=0$, the factor $(z/i)^\chi$ is interpreted as identically one.
This representation also yields Theorem~\ref{thm_2}, by convergence
of finite products and a uniform estimate for the tails.
For $|z|\le r$,
\[
 \frac{1-z^2/t_{a,n}^2}{1+1/t_{a,n}^2}
       =1-\frac{z^2+1}{t_{a,n}^2+1},\qquad
 \sum_{n>N}\frac{|z^2+1|}{t_{a,n}^2+1}
       \le(r^2+1)\sum_{n>N}\gamma_n^{-2}\longrightarrow0
\]
uniformly in $a$.
The inequality
$|\prod_k(1+b_k)-1|\le\exp(\sum_k|b_k|)-1$
then shows that the product tails tend uniformly to one.
This is the approach used in the proof of
\cite[Theorem 1.2]{Su26b}.

\subsection{A normalized limit involving the function field zeta function}
\label{sec:finite-corollary}

First, put
\begin{equation} \label{eq_616}
 G:=A/\Psf,\quad E:=A+iA',\quad
 E_0:=E/G=\Psf+iB_0,\quad B_0:=A'/G.
\end{equation}
The function $G$ is entire and has a zero of order $m_\gamma-1$ at
each $\gamma \in \Gamma$.
Thus $B_0$ is also entire, and
\begin{equation} \label{eq_617}
 B_0(\gamma)=m_\gamma\Psf'(\gamma)\ne0.
\end{equation}
Since $\Psf$ and $B_0$ are real-valued on the real axis and have no
common real zeros, $E_0$ has no real zeros. Moreover,
\begin{equation} \label{eq_618}
 -\frac{A'(z)}{A(z)}
 =-\frac h2\sum_{j=1}^dm_j\cot\frac{h(z-\beta_j)}2
\end{equation}
has positive imaginary part in the upper half-plane.
Using $E^\sharp/E=(1-iA'/A)/(1+iA'/A)$, we obtain
$|E^\sharp(z)|<|E(z)|$ for $z\in\C_+$.
Hence $E_0$ is a Hermite--Biehler function without real zeros.

The generating function $E_\ast$ of $\widehat{\Hstar}$ constructed
above generally differs from $E=A+iA'$ and $E_0=E/G$.
In particular, $B_\ast$ and $B_0$ generally differ.
The existence of self-adjoint extensions with the same zero set therefore
does not identify the norms of the two spaces or their full families of
self-adjoint extensions.

\begin{theorem} \label{thm_607}
For $W(a,\theta;z)$ defined by \eqref{eq_510} and the phase in \eqref{eq_105}, 
we have 
\begin{equation} \label{eq_619}
 \lim_{a\uparrow a_\ast}
 \frac{\pi}{2B_\ast(i)E_0(z)}W(a,\theta_\ast;z)
 =\frac{A(z)}{E(z)}, 
\end{equation}
where the convergence is uniform on every compact subset of $\C_+$.
\end{theorem}
\begin{proof}
By Theorem~\ref{thm_602} and the identities
$A=G\Psf$, $E=GE_0$, the left-hand side converges to
$\Psf(z)/E_0(z)=A(z)/E(z)$.
Since $B_\ast(i)\ne0$ and $E_0$ has no zeros in the upper half-plane,
the convergence is uniform on every compact subset of that half-plane.
\end{proof}

Since $E_0$ has no zeros in $\C_+$, we can choose
a holomorphic logarithm of $E_0$ there. Thus
\begin{equation} \label{eq_620}
 \phi(z):=\log\!\left(\frac{\pi}{2B_\ast(i)}\right)-\log E_0(z)
\end{equation}
is holomorphic in $\C_+$, and \eqref{eq_619}
becomes $e^{\phi(z)}W(a,\theta_\ast;z)\to A(z)/E(z)$.  
For a curve, put $E_C:=A_C+iA_C'$.
With $\xi_C'$ denoting differentiation with respect to $s$, we have
\[
 \frac{A_C(z)}{E_C(z)}
 =\frac{\xi_C(1/2-iz)}
        {\xi_C(1/2-iz)+\xi_C'(1/2-iz)}.
\]
This is the function field analogue of the limit in
\cite[Corollary~1.6]{Su26a}.
If $\Psf$ is replaced by $c\Psf$ with $c\in\R\setminus\{0\}$,
then $E_0$ is replaced by $cE_0$ and $B_\ast$ by $B_\ast/c$,
so the factor in \eqref{eq_619} is unchanged.

\begin{corollary} \label{cor_608}
We have
\begin{equation} \label{eq_621}
 \lim_{a\uparrow a_\ast}
 \frac{\Psf(i)}{E_0(z)}
 \frac{W(a,\theta_\ast;z)}{W(a,\theta_\ast;i)}
 =\frac{A(z)}{E(z)},
\end{equation}
where the convergence is uniform on every compact subset of $\C_+$.
In the case of a curve, if all zeros of $A_C$ are simple, then
\begin{equation} \label{eq_622}
 \lim_{a\uparrow a_\ast} \frac{W(a,\pi;z)}{W(a,\pi;i)} =\frac{A_C(z)}{A_C(i)},
\end{equation}
where the convergence is uniform on every compact subset of $\C$.
In this case, $\Psf=A_C/A_C(0)$, $G=A_C(0)$, and
$E_0=E_C/A_C(0)$.
\end{corollary}
\begin{proof}
Multiply the limit in Theorem~\ref{thm_2} by $\Psf(i)/E_0(z)$
to obtain the first formula.
The second assertion is the part of that theorem concerning simple zeros.
\end{proof}

If $0\notin\Gamma$, then $W(a,\pi;0)\ne0$, so normalization at the
origin is possible.
If $0\in\Gamma$, however, $W(a,0;0)=0$.
We have therefore used $i$ as the normalization point in
\eqref{eq_107}, since it is available in both cases.

\section{Global realizations} \label{section_7}

In the proofs of Theorems~\ref{thm_1} and \ref{thm_2}, the space
$\mathcal H(A_{a_\ast})$ has played the role of the global Hilbert
space $\Hglob$ in \cite{Su26a, Su26b}.
This suffices for the main results of the present paper. 
Nevertheless, a counterpart of $\Hglob$ can also be defined for
function fields.
We describe the relation between the two spaces to clarify how the
theory differs from that for the Riemann zeta function.
The resulting global realization also provides the basis for the
geometric and operator-theoretic interpretations developed in
Sections~\ref{section_8} and \ref{section_9}.

\subsection{The global Hilbert space}

Let $N_\infty:=\{u\in C_c^\infty(\R):Q_C(u)=0\}$ and define
\[
 \Hglob:=\overline{C_c^\infty(\R)/N_\infty}^{\,Q_C}.
\]
By \eqref{eq_104}, the map sending the class $[u]$ of a test function
to $(\widehat u(\gamma))_{\gamma\in\Gamma}$ is an isometry into
$\ell^2(\Gamma,m)$.
Theorem~\ref{thm_401} and the density of test functions on the
endpoint interval imply that its range is dense in $\ell^2(\Gamma,m)$.
It therefore extends uniquely by completion to an isometric isomorphism
\begin{equation} \label{eq_701}
 \Samp:\Hglob\overset{\sim}{\longrightarrow}\ell^2(\Gamma,m),
 \qquad
 \Samp[u]=(\widehat u(\gamma))_{\gamma\in\Gamma}.
\end{equation}

Let $j_a:\Ha\to\Hglob$ be the natural isometric embedding, and put
$\mathcal M_a=j_a\Ha$ and $\Pi_a=j_aj_a^\ast$.
Then $\Pi_a$ is the orthogonal projection onto $\mathcal M_a$.
The map $j_{a_\ast}$ is surjective, and $\Pi_a\to I$ strongly as
$a\uparrow a_\ast$.
For $a>a_\ast$, the map $j_a$ defined through the quotient and
completion is also surjective.
Thus every global element has a unique $L^2$ representative supported
in the interval at the endpoint, whereas uniqueness of such
representatives fails on larger intervals.
\medskip

On $\ell^2(\Gamma,m)$, define
\[
 (T_\Gamma c)_\gamma=\gamma c_\gamma,\qquad
 \Dom(T_\Gamma)=\left\{c\in\ell^2(\Gamma,m):
                         \sum_\gamma m_\gamma\gamma^2
                         |c_\gamma|^2<\infty\right\},
\]
and set
\begin{equation} \label{eq_702}
 \Dlim:=\Samp^{-1}T_\Gamma\Samp.
\end{equation}
This is self-adjoint, its spectrum is $\Gamma$, and every eigenvalue
is simple.
For test functions, $\Dlim[u]=[iu']$.
Indeed, $\sum m_\gamma\gamma^2|\widehat u(\gamma)|^2<\infty$,
and $N_\infty$ is invariant under differentiation.
At $a=a_\ast$, Lemma~\ref{lem_603}
and the endpoint sampling isomorphism give 
$
 \Dlim=j_{a_\ast}
 \overline\D_{a_\ast,\theta_\ast}
 j_{a_\ast}^{-1}$, including the domains.

\subsection{de Branges spaces and the model space}

With $G$ and $E_0$ as in \eqref{eq_616}, we have
\[
 \mathcal H(E)=G\mathcal H(E_0),\qquad
 \|GF\|_{\mathcal H(E)}=\|F\|_{\mathcal H(E_0)}.
\]
Set $\Theta:=E^\sharp/E=E_0^\sharp/E_0$.
The corresponding model space is
\[
 \mathcal K:=K(\Theta)
 =H^2(\C_+)\ominus\Theta H^2(\C_+)
 =\mathcal H(E)/E,
\]
with the usual boundary $L^2(\R)$ norm.
The map $T_EF:=EF$ is an isometric isomorphism from $\mathcal K$
to $\mathcal H(E)$.
For basic facts about model spaces, see \cite[Section~3.1]{Su23} and
\cite[Sections~2.3--2.4]{Su25}.

The identity $A/E=(1+\Theta)/2$ and Taylor expansion give
$\Theta(\gamma)=-1$ and $\Theta'(\gamma)=-2i/m_\gamma$.
Define
\begin{equation} \label{eq_703}
 F_\gamma(z):=i\sqrt{\frac{m_\gamma}{\pi}}
                \frac{A(z)}{E(z)(z-\gamma)},\qquad
 F_\gamma(\eta)=\frac{\delta_{\gamma\eta}}{\sqrt{\pi m_\gamma}}.
\end{equation}
Then $\{F_\gamma\}_{\gamma\in\Gamma}$ is an orthonormal basis of
$\mathcal K$.
Indeed, by de Branges' basis theorem \cite[Theorem~22]{dB68}, the only
possible missing direction occurs when $\Psf\in\mathcal H(E_0)$.
Here $\Psf/E_0=A/E$ is a nonzero continuous periodic boundary
function, so it does not belong to $L^2(\R)$, and no direction is
missing.
The squared norm of each reproducing kernel is
$|\Theta'(\gamma)|/(2\pi)=1/(\pi m_\gamma)$, giving the normalization
in \eqref{eq_703}.
This agrees with \cite[Proposition~3.2]{Su23}; completeness has been
verified here for the present abstract periodic data.
Moreover,
\begin{equation} \label{eq_704}
 U:\Hglob\overset{\sim}{\longrightarrow}\mathcal K,\qquad
 U\xi:=\sum_\gamma\sqrt{m_\gamma}(\Samp \xi)_\gamma F_\gamma,
 \qquad U_a:=Uj_a
\end{equation}
is unitary.
Let $\widehat{f_\gamma}=F_\gamma$ and set
\[
V:=\mathcal F^{-1}\mathcal K\subset L^2(0,\infty).
\]
Through this isomorphism, associate with each $\xi\in\Hglob$ the
function $v_\xi\in V$ defined by
\begin{equation} \label{eq_705}
 v_\xi:=\sqrt\pi\,\mathcal F^{-1}U\xi
     =\sqrt\pi\sum_\gamma
       \sqrt{m_\gamma}(\Samp \xi)_\gamma f_\gamma,
 \qquad
 2\|v_\xi\|_2^2=\|\xi\|_{\Hglob}^2.
\end{equation}
Then $\widehat{v_\xi}(\gamma)=(\Samp \xi)_\gamma$, and
\[
 \|E\widehat{v_\xi}\|_{\mathcal H(E)}^2
       =\|\widehat{v_\xi}\|_2^2
       =2\pi\|v_\xi\|_2^2
       =\pi\|\xi\|_{\Hglob}^2.
\]
This corresponds to the relation in \cite[Theorem~5.5]{Su25}.
Point values are understood through analytic continuation of functions
in the model space, not by assigning pointwise values to the Fourier
integrals of arbitrary $L^2$ functions.
Although the compactly supported representative at the endpoint and
$v_\xi$ have the same sampling values, their Fourier transforms need not
agree at general complex points.

Define multiplication by the independent variable on $\mathcal H(E)$
by
$\Dom(M):=\{F\in\mathcal H(E):zF(z)\in\mathcal H(E)\}$
and $(MF)(z):=zF(z)$.
In $\mathcal H(E_0)$, where the real zeros have been removed, write
$s_\alpha$ for the associated function.
For every $\alpha$, the function
$s_\alpha/E_0=\frac i2(e^{i\alpha}-e^{-i\alpha}\Theta)$
is a nonzero periodic boundary function and does not belong to $L^2$.
Hence $M$ is densely defined; see \cite[Section~2.3]{Su25} for this
criterion and the description of its extensions.
Let $M_{\pi/2}$ be the extension whose spectrum is the zero set of
$\Psf$.
Then
\[
 M_{\pi/2}(EF_\gamma)=\gamma EF_\gamma,\qquad
 (T_EU)\Dlim(T_EU)^{-1}=M_{\pi/2}.
\]
The domains also correspond under the eigenvector expansions.

\subsection{Relation to compression at the finite endpoint}

We describe the relation between the formulation of
Section~\ref{section_6-compression} and the global formulation.
Surjectivity of sampling at the endpoint implies that
$j_{a_\ast}:\Hstar\to\Hglob$ is an isometric isomorphism, and
\begin{equation} \label{eq_706}
 j_a=j_{a_\ast}\iota_a,\qquad
 \Dlim=j_{a_\ast}D_\ast j_{a_\ast}^{-1},
\end{equation}
where $\iota_a$ is the extension by zero used in
Section~\ref{section_6-compression}.

\begin{theorem} \label{thm_701}
If $0\notin\Gamma$, then, for every $0<a\le a_\ast$,
\[
 (\overline\D_{a,\pi})^{-1}=j_a^\ast\Dlim^{-1}j_a
\]
with equality of domains.
\end{theorem}
\begin{proof}
By \eqref{eq_706} and
Proposition~\ref{prop_604}, the right-hand side equals
$\iota_a^\ast D_\ast^{-1}\iota_a$.
\end{proof}

If $0\in\Gamma$, put $\Hred:=\Hglob\ominus\ker\Dlim$.
The map $j_{a_\ast}|_{V_\ast}:V_\ast\to\Hred$ is unitary, and
$j_{a,V}:=j_{a_\ast}\kappa_a:V_a\to\Hred$ is an isometric embedding.
Set
$\mathsf C:=(j_{a_\ast}|_{V_\ast})\mathsf C_\ast
(j_{a_\ast}|_{V_\ast})^{-1}$.

\begin{theorem} \label{thm_702}
If $0\in\Gamma$, then, for every $0<a\le a_\ast$,
\[
 (D_a|_{V_a})^{-1}=j_{a,V}^\ast\mathsf Cj_{a,V}
\]
with equality of domains.
\end{theorem}

\begin{proof}
Conjugate Proposition~\ref{prop_605} by the
unitary map $j_{a_\ast}|_{V_\ast}$.
\end{proof}

These global formulas are not used in the proofs of the main theorems.
They describe the correspondence between the construction at the finite
endpoint $a=a_\ast$ and the construction for the Riemann zeta function
in \cite{Su26a, Su26b}.

\subsection{Convergence of compressed inverses and resolvents}

Let $P_a^V$ be the orthogonal projection in $V_\ast$ onto
$\kappa_a V_a$.
Then $P_a^V\to I$ strongly.
If there is a central zero, density follows by approximating a mean-zero
function at the endpoint by test functions and correcting their integrals
with a fixed test function of integral one.
By \eqref{eq_403}, the same approximation holds in the
$Q_C$-norm.
Put
\[
 \widetilde{\mathsf C}_a:=P_a^V\mathsf C_\ast P_a^V
 =\kappa_a\mathsf C_a\kappa_a^\ast.
\]
Since $\mathsf C_\ast$ is Hilbert--Schmidt,
\begin{equation} \label{eq_707}
 \|\widetilde{\mathsf C}_a-\mathsf C_\ast\|_{\mathfrak S_2}
 \longrightarrow0
 \qquad(a\uparrow a_\ast).
\end{equation}
Indeed, for an orthonormal basis $(e_n)$, dominated convergence gives
\[
 \sum_n\|(I-P_a^V)\mathsf C_\ast e_n\|^2\to0.
\]
Taking adjoints and using the triangle inequality proves the assertion.
This is a stronger form of operator convergence that was not needed in
Section~\ref{sec:finite-minmax}.

\begin{corollary} \label{cor_703}
Let $K\subset\C\setminus\Gamma$ be compact.
For $a<a_\ast$ sufficiently close to $a_\ast$, we have
$K\subset\rho(D_a)$, and
$\iota_a(D_a-z)^{-1}\iota_a^\ast~\to~(D_\ast-z)^{-1}$
in operator norm, uniformly for $z\in K$.
By \eqref{eq_706}, on the global space we likewise have
$j_a(D_a-z)^{-1}j_a^\ast~\to~(\Dlim-z)^{-1}$
in operator norm, uniformly for $z\in K$.
\end{corollary}
\begin{proof}
If $0\notin\Gamma$, apply \eqref{eq_707} and stability of
inversion to
$\widetilde{\mathsf C}_a(I-z\widetilde{\mathsf C}_a)^{-1}$.
If $0\in\Gamma$, extend $\widetilde{\mathsf C}_a$ by zero to
$\Hstar$.
Let $q_a$ be the orthogonal projection onto the span of
$\iota_au_{0,a}/\|u_{0,a}\|_{A_a}$.
The compressed resolvent is then
$\widetilde{\mathsf C}_a(I-z\widetilde{\mathsf C}_a)^{-1}
-z^{-1}q_a$.
By \eqref{eq_605}, $q_a$ converges in operator norm to the
projection onto the zero eigenspace at the endpoint.
Since $0\notin K$, the conclusion follows in the same way.
\end{proof}

\section{A geometric proof of positivity of the Weil quadratic form} \label{section_8}

We describe the relation between the quadratic form $Q_C$ and positivity of
the Rosati involution on the Jacobian of the curve.
The results of this section are not used in proving
Theorems~\ref{thm_1} and \ref{thm_2}, but the Hilbert space
realization obtained here is used in Section~\ref{section_9}.
Following the constructions of Swinnerton-Dyer~\cite{SD71},
Patterson~\cite{Pa88}, and Tretkoff~\cite{Tr06},
we express $Q_C$ in terms of the Rosati trace pairing.

Let $J=J(C)$ be the Jacobian of $C$, and let $\dagger$ denote
the Rosati involution associated with its canonical principal
polarization.
Put $B=\operatorname{End}^0_{\mathbb F_q}(J):=
\operatorname{End}_{\mathbb F_q}(J)\otimes_{\Z}\mathbb Q$.
For a prime $\ell\ne\operatorname{char}\mathbb F_q$, define
$
V_\ell:=
\bigl(\varprojlim_n J[\ell^n](\overline{\mathbb F}_q)\bigr)
\otimes_{\mathbb Z_\ell}\mathbb Q_\ell
$,
where $J[\ell^n]$ is the $\ell^n$-torsion subgroup of $J$.
Define the trace on $B$ by
\[
 \tau(x):=\operatorname{Tr}\bigl(x\mid V_\ell\bigr)
 \qquad(x\in B).
\]
It is rational-valued and independent of $\ell$, and extends
linearly to $B_{\R}:=B\otimes_{\mathbb Q}\R$ and
$B_{\C}:=B\otimes_{\mathbb Q}\C$.
On $B_{\C}$, write $\star$ for the conjugate-linear extension
of $\dagger$:
\[
 (x\otimes z)^\star:=x^\dagger\otimes\overline z
 \qquad(x\in B,\ z\in\C).
\]
Positivity of the Rosati involution gives
\begin{equation} \label{eq_801}
 \tau(XX^\star)\ge0\qquad(X\in B_{\C}),
\end{equation}
with equality only when $X=0$
(\cite[Theorem~9 and p.~46]{SD71},
\cite[Theorem~17.3]{Mi86}).
Indeed, writing $X=X_1+iX_2$ with $X_1,X_2\in B_{\R}$,
symmetry of the Rosati trace pairing gives
$\tau(XX^\star)=
\tau(X_1X_1^\dagger)+\tau(X_2X_2^\dagger)$.

We recall the relation of this positivity to intersection theory.
Put $\bar C=C\otimes_{\mathbb F_q}\overline{\mathbb F}_q$,
and denote the transpose of a correspondence $D$ on
$\bar C\times\bar C$ by $D^t$.
Let $\Delta$ be the diagonal correspondence and define
$d(D):=D\cdot(\{P\}\times\bar C)$, where $P\in\bar C$ is any
point and $\cdot$ denotes the intersection number on
$\bar C\times\bar C$.
Swinnerton-Dyer~\cite[Section~7]{SD71} constructs a surjection
\[
 \Phi:\operatorname{Corr}(\bar C)
       \longrightarrow
       \operatorname{End}_{\overline{\mathbb F}_q}(\bar J).
\]
Here $\operatorname{Corr}(\bar C)$ is the ring of
correspondences given by rational equivalence classes of
divisors, and
$\bar J=J\otimes_{\mathbb F_q}\overline{\mathbb F}_q$.
Its kernel consists of classes of horizontal and vertical
correspondences (\cite[(27) and Lemma~2]{SD71}).
He also defines
\begin{equation} \label{eq_802}
 S(D):=d(D)+d(D^t)-D\cdot\Delta.
\end{equation}
Under the duality induced by the canonical principal polarization
of the Jacobian, transposition of correspondences corresponds
to the Rosati involution.
By the Lefschetz trace formula for correspondences, $S(D)$ equals
the trace of $\Phi(D)$ on the Tate module.
Consequently,
\begin{equation} \label{eq_803}
 \begin{gathered}
 \Phi(D^t)=\Phi(D)^\dagger,\qquad
 S(D)=\tau(\Phi(D)),\qquad
 S(D\circ D^t)
   =\tau\bigl(\Phi(D)\Phi(D)^\dagger\bigr), 
 \end{gathered}
\end{equation}
where $\tau$ is extended to the geometric endomorphism algebra
by taking the trace on the same Tate module.
For the diagonal correspondence, $S(\Delta)=2g=\tau(1)$.
In particular, \cite[Theorem~9]{SD71} proves that 
$S(D\circ D^t)>0$ when $\Phi(D)\ne0$ by intersection theory.
Thus \eqref{eq_803} expresses the relation
between that geometric positivity and the Rosati trace pairing.

Let $\pi\in\operatorname{End}_{\mathbb F_q}(J)$ be the
$q$-Frobenius endomorphism, and put
\[
 F:=q^{-1/2}\pi\in B_{\R}.
\]
The identity $\pi^\dagger\pi=[q]$ gives $F^\dagger=F^{-1}$
(\cite[Lemma~19.2]{Mi86},
\cite[Section~5.17]{Pa88}).

For a curve, the numbers $p_n$ are determined by
\eqref{eq_302}, together with $p_0=2g$ and
$p_{-n}=p_n$ for $n\ge1$.
The Lefschetz trace formula for Frobenius gives
$\tau(\pi^n)=q^n+1-N_n$ for $n\ge1$ and $\tau(1)=2g$.
Moreover, $\tau(Y^\dagger)=\tau(Y)$ and
$F^\dagger=F^{-1}$ imply $\tau(F^{-n})=\tau(F^n)$.
Therefore,
\begin{equation} \label{eq_804}
 \tau(F^n)=p_n\qquad(n\in\Z).
\end{equation}

For a finitely supported function $a:\Z\to\C$, put
\[
 a^\vee(k):=\overline{a(-k)},\qquad
 M(a,z):=\sum_{k\in\Z}a(k)z^k,\qquad
 X_a:=\sum_{k\in\Z}a(k)F^k\in B_{\C}.
\]
Write $P_C(T)=\prod_{\nu=1}^{2g}(1-\omega_\nu T)$,
counting the $\omega_\nu$ with multiplicity, and put
$\alpha_\nu=\omega_\nu/\sqrt q$.
With the arrangement of zeros in \eqref{eq_202},
our $p_n$ equals $\sum_\nu\alpha_\nu^{-n}$.
The symmetry of reciprocal roots implied by the functional
equation gives
\[
 \sum_{\nu=1}^{2g}\alpha_\nu^n
 =\sum_{\nu=1}^{2g}\alpha_\nu^{-n}
 =p_n\qquad(n\in\Z).
\]
For finitely supported functions $a,b:\Z\to\C$, the
Hermitian form of Tretkoff~\cite[(23)]{Tr06} is
\[
 R(a,b):=\sum_{\nu=1}^{2g}M(a\ast b^\vee,\alpha_\nu),
\]
where $\ast$ denotes discrete convolution on $\Z$.
Expanding the finite sums and using \eqref{eq_804},
we obtain
\begin{equation} \label{eq_805}
 \begin{split}
 R(a,b)
  =\sum_{k,\ell\in\Z}
      a(k)\overline{b(\ell)}\,p_{k-\ell}
  =\tau(X_aX_b^\star).
 \end{split}
\end{equation}
In particular, \eqref{eq_801} gives $R(a,a)\ge0$.
This expresses Patterson's construction from finitely supported
functions and Laurent polynomials in normalized Frobenius in
terms of Tretkoff's discrete quadratic form.

\begin{proposition} \label{prop_801}
Let $u,v\in C_c^\infty(\R)$. For $0\le t<h=\log q$,
we define
\[
 u_t(k):=u(t+kh),\qquad
 v_t(k):=v(t+kh),\qquad
 X_u(t):=\sum_{k\in\Z}u(t+kh)F^k.
\]
Then
\begin{equation} \label{eq_806}
 Q_C(u,v)
 =h\int_0^hR(u_t,v_t)\,dt
 =h\int_0^h\tau\bigl(X_u(t)X_v(t)^\star\bigr)\,dt.
\end{equation}
Consequently, $Q_C(u)\geq0$.
\end{proposition}

\begin{proof}
For each $t\in[0,h)$, the functions $u_t$ and $v_t$ have
finite support, and these supports lie in a fixed finite set
independent of $t$.
We may therefore integrate \eqref{eq_805}
term by term:
\begin{align*}
 h\int_0^h R(u_t,v_t)\,dt
 &=h\sum_{k,\ell\in\Z}p_{k-\ell}
      \int_0^h u(t+kh)\overline{v(t+\ell h)}\,dt\\
 &=h\sum_{n\in\Z}p_n
      \int_{\R}u(x)\overline{v(x-nh)}\,dx=h\sum_{n\in\Z}p_n(u*\widetilde v)(nh),
\end{align*}
where $v(x):=\overline{v(-x)}$. 
The last sum is also finite.
Equation~\eqref{eq_301} and 
the explicit formula 
$Q_C(u,v)=\Wdist(u*\widetilde v)$ give the first equality.
The second follows from \eqref{eq_805}, and
positivity follows from \eqref{eq_801}.
\end{proof}

This derivation of positivity does not use the fact that the
zeros lie on the critical line.
Indeed, if one first defines
$h\sum_np_n(u*\widetilde v)(nh)$ using the $p_n$ in
\eqref{eq_302}, its positivity follows solely from
\eqref{eq_801} and \eqref{eq_804}.

The identity \eqref{eq_806} also admits an interpretation
in terms of intersection theory.
Extend the correspondence ring in
\eqref{eq_803} to rational coefficients,
divide by the kernel of $\Phi$, and then extend scalars to
$\R$ and $\C$.
Let $\Pi$ be the class of the graph of Frobenius and put
$\mathfrak F:=q^{-1/2}\Pi$.
Then $\Phi(\mathfrak F)=F$ and
$\mathfrak F^t=\mathfrak F^{-1}$, where the transpose and
inverse are taken in the quotient ring by $\ker\Phi$.
On this quotient, also write $\star$ for the conjugate-linear
extension of transposition.
Define
\[
 D_u(t):=\sum_{k\in\Z}u(t+kh)\mathfrak F^k.
\]
Since $\Phi(D_u(t))=X_u(t)$, equations
\eqref{eq_803} and \eqref{eq_806} give
\begin{equation} \label{eq_807}
 Q_C(u,v)
 =h\int_0^h S\bigl(D_u(t)\circ D_v(t)^\star\bigr)\,dt, 
\end{equation}
where $S$ is extended linearly to correspondence classes and
to the enlarged coefficient fields.
Thus the positivity $Q_C(u)\ge0$ can also be viewed as an application of
the positivity from intersection theory in
Swinnerton-Dyer~\cite[Theorem~9]{SD71}.

\begin{corollary} \label{cor_802}
Let $\mathcal E:=\operatorname{span}_{\C}\{F^k:k\in\Z\}$
with the inner product
\[
 \langle X,Y\rangle_{\mathcal E}:=\tau(XY^\star),
 \qquad X,Y\in\mathcal E.
\]
The induced inner product on $L^2((0,h),\mathcal E)$ is 
\[
 \langle f,g\rangle_{L^2}
 :=\int_0^h\langle f(t),g(t)\rangle_{\mathcal E}\,dt.
\]
Then the map
\[
 C_c^\infty(\R)~\to~L^2((0,h),\mathcal E),
 \qquad
 u\mapsto
 \bigl(t\mapsto\sqrt h\,X_u(t)\bigr)
\]
passes to the quotient by the null space of $Q_C$ and extends
by completion to an isometric isomorphism
\begin{equation} \label{eq_808}
 \mathcal U:\Hglob\overset{\sim}{\longrightarrow}L^2((0,h),\mathcal E).
\end{equation}
\end{corollary}
\begin{proof}
The isometry follows from \eqref{eq_806}. 
Choose $\psi\in C_c^\infty(0,h)$ and
$X=\sum_{k\in I}c_kF^k\in\mathcal E$, where $I$ is finite.
Set $u(t+kh)=c_k\psi(t)$ for $0<t<h$ and $k\in I$,
and set $u=0$ elsewhere.
Then $u\in C_c^\infty(\R)$ and $X_u(t)=\psi(t)X$.
Linear combinations of these $\mathcal E$-valued functions
are dense in $L^2((0,h),\mathcal E)$, so the isometric
embedding given by \eqref{eq_806} extends to
a surjection.
\end{proof}

\section{The differentiation operator and Frobenius} \label{section_9}

We now relate the differentiation operator at the endpoint to Frobenius.
The essential intermediate object is the strongly continuous unitary
group of translations on the global Hilbert space.
For $s\in\R$ and $u\in C_c^\infty(\R)$, set
$(T_su)(x)=u(x-s)$ and $\mathsf T_s[u]=[T_su]$.

\begin{proposition} \label{prop_901}
The operators $\mathsf T_s$ extend to a strongly continuous
one-parameter unitary group on $\Hglob$, and
\begin{equation} \label{eq_901}
 \mathsf T_s=e^{is\Dlim}\qquad(s\in\R),
\end{equation}
where $\Dlim$ is the self-adjoint operator defined in \eqref{eq_702}.
\end{proposition}

\begin{proof}
We observe that translations preserve the convolution
$u*\widetilde v$, hence $Q_C(T_su,T_sv)=Q_C(u,v)$. 
Thus they preserve the null space of $Q_C$ and induce
unitary operators on its completion.
Moreover,
$\widehat{T_su}(\gamma)=e^{is\gamma}\widehat u(\gamma)$. 
For $c=\Samp\xi$, we obtain from the sampling isomorphism
\[
 \|\mathsf T_s\xi-\xi\|_{\Hglob}^2
 =\sum_{\gamma\in\Gamma}m_\gamma
 |e^{is\gamma}-1|^2|c_\gamma|^2~\to~0
\]
by dominated convergence.
Thus the group is strongly continuous. 
Under $\Samp$, the translation group acts diagonally
by multiplication by $e^{is\gamma}$, giving \eqref{eq_901}.
We see that its self-adjoint generator has domain
\[
 \left\{c\in\ell^2(\Gamma,m):
   \sum_{\gamma\in\Gamma}m_\gamma\gamma^2|c_\gamma|^2<\infty
 \right\}
\]
in sampling coordinates, which is exactly the domain of $\Dlim$.
Thus $\mathsf T_s=e^{is\Dlim}$, and its infinitesimal
generator is $i\Dlim$.
\end{proof}

Suppose now that $C/\mathbb F_q$ is a curve.
We use the unitary map
$\mathcal U:\Hglob\to L^2((0,h),\mathcal E)$
defined in \eqref{eq_808}, where
$F=q^{-1/2}\pi$ and
$\mathcal E=\operatorname{span}_{\C}\{F^k:k\in\Z\}$.
We define the operator $M_F$ on
$L^2((0,h),\mathcal E)$ by
\[
 (M_Ff)(t)=Ff(t).
\]
\begin{proposition} \label{prop_902}
The operator $M_F$ is unitary and
\begin{equation} \label{eq_902}
 \mathcal Ue^{ih\Dlim}\mathcal U^{-1}=M_F.
\end{equation}
More precisely, for $f\in L^2((0,h),\mathcal E)$,
we define its extension to $\R$ by
\[
 \widetilde f(t+nh)=F^{-n}f(t),
 \quad 0<t<h,\quad n\in\Z.
\]
Then
\begin{equation} \label{eq_903}
 (\mathcal U\mathsf T_s\mathcal U^{-1}f)(t)
 =\widetilde f(t-s)
 \quad\text{a.e. }t\in(0,h).
\end{equation}
For $0<s<h$ this equals $f(t-s)$ when $s<t<h$ and
$Ff(t-s+h)$ when $0<t<s$.
\end{proposition}

\begin{proof}
Since $F^\star=F^{-1}$ and the trace is cyclic,
\[
 \langle FX,FY\rangle_{\mathcal E}
 =\tau(FXY^\star F^{-1})
 =\tau(XY^\star)
 =\langle X,Y\rangle_{\mathcal E}.
\]
Thus $M_F$ is unitary.
For test functions $u \in C_c^\infty(\R)$ the sum
$X_u(t)=\sum_{k\in\Z}u(t+kh)F^k$ is locally finite and
satisfies
\[
 X_u(t+h)=F^{-1}X_u(t),\qquad
 X_{T_su}(t)=X_u(t-s).
\]
This proves \eqref{eq_903} on a dense subspace.
Its right-hand side defines a unitary operator on the entire
$L^2$ space, since unitarity of $F$ makes
$\|\widetilde f(t)\|_{\mathcal E}^2$ periodic of period $h$.
By density and continuity, the identity holds for
every $f\in L^2((0,h),\mathcal E)$. 
For $s=h$, the twisted periodicity gives
$\widetilde f(t-h)=Ff(t)$.
Combining this with \eqref{eq_901},
we obtain \eqref{eq_902}.
\end{proof}

\begin{theorem} \label{thm_903} 
We define
\[
 D_F f:=if',\qquad 
 \Dom(D_F)=
 \{f\in H^1((0,h),\mathcal E):f(h)=F^{-1}f(0)\}.
\]
Then $D_F$ is self-adjoint and
\begin{equation} \label{eq_904}
 \mathcal U\Dlim\mathcal U^{-1}=D_F
\end{equation}
as equality of operators including their domains.
In particular, if
$D_\ast=\overline\D_{a_\ast,\theta_\ast}$, then
\begin{equation} \label{eq_905}
 (\mathcal Uj_{a_\ast})D_\ast
  (\mathcal Uj_{a_\ast})^{-1}=D_F,
 \qquad
 (\mathcal Uj_{a_\ast})e^{ihD_\ast}
  (\mathcal Uj_{a_\ast})^{-1}=M_F.
\end{equation}
\end{theorem}
\begin{proof}
Since $\mathcal E$ is a finite-dimensional inner-product space, 
the standard results for vector-valued Sobolev spaces apply to $H^1((0,h),\mathcal E)$. 
For $f,g\in H^1((0,h),\mathcal E)$, integration by parts gives
\[
 \langle if',g\rangle_{L^2}-\langle f,ig'\rangle_{L^2}
 =i\{\langle f(h),g(h)\rangle_{\mathcal E}
       -\langle f(0),g(0)\rangle_{\mathcal E}\}.
\]
Since $F$ is unitary, the right-hand side vanishes for $f,g\in\Dom(D_F)$. 
Conversely, we see that
every $f\in\Dom(D_F^*)$ has a weak derivative in
$L^2((0,h),\mathcal E)$ by testing the adjoint against
$C_c^\infty((0,h),\mathcal E)$, and hence belongs to
$H^1((0,h),\mathcal E)$.
We then test against functions with arbitrary value $v$ at zero
and value $F^{-1}v$ at $h$ to obtain
$f(h)=F^{-1}f(0)$.
Thus $D_F$ is self-adjoint. 

For $f\in\Dom(D_F)$, we see that the endpoint traces
of the extension $\widetilde f$ in \eqref{eq_903}
agree at every multiple of $h$.
Hence $\widetilde f$ belongs to
$H^1_{\mathrm{loc}}(\R,\mathcal E)$.
We therefore obtain
\[
 \frac{\widetilde f(\,\cdot-s)-f}{s}
 \longrightarrow -f'=iD_F f
 \quad\text{in }L^2((0,h),\mathcal E).
\]
In view of \eqref{eq_901} and \eqref{eq_903},
we conclude that
$D_F\subset\mathcal U\Dlim\mathcal U^{-1}$.
Both sides are self-adjoint, so equality follows,
including equality of their domains.
Finally, we obtain the first equality in \eqref{eq_905}
from \eqref{eq_706} and \eqref{eq_904},
and the second from \eqref{eq_902},
since unitary equivalence of self-adjoint operators
is preserved under functional calculus.
\end{proof}

\begin{lemma} \label{lem_904}
There is an isometric algebra isomorphism
\[
 \mathcal E~\to~\C^d,\qquad
 P(F)\mapsto(P(\zeta_1),\ldots,P(\zeta_d)),
\]
where $\zeta_j=e^{ih\beta_j}$, and the right-hand side
has inner product
$\langle x,y\rangle=\sum_{j=1}^d m_jx_j\overline{y_j}$.
Consequently $\dim_{\C}\mathcal E=d$, and multiplication
by $F$ on $\mathcal E$ has simple eigenvalues $\zeta_j$.
\end{lemma}
\begin{proof}
By \eqref{eq_804} and the periodic data, we have
\[
 \tau(P(F)Q(F)^\star)
 =\sum_{j=1}^dm_jP(\zeta_j)\overline{Q(\zeta_j)}
\]
for Laurent polynomials $P,Q$.
Since the Rosati trace pairing is strictly positive,
the evaluation map is well defined and injective.
By Lagrange interpolation at the distinct points
$\zeta_1,\ldots,\zeta_d$, we see that it is also surjective.
Thus it is an isometric algebra isomorphism.
Under this isomorphism, multiplication by $F$
acts as multiplication by $\zeta_j$ in the $j$th coordinate.
\end{proof}

Let $b_j\in\mathcal E$ correspond to the $j$th standard
coordinate vector under the isomorphism in Lemma \ref{lem_904}.
Then $Fb_j=\zeta_jb_j$ and $\|b_j\|_{\mathcal E}^2=m_j$.
For $1\le j\le d$ and $n\in\Z$, put
\[
 \phi_{j,n}(t)
 :=\frac{e^{-i(\beta_j+2\pi n/h)t}}{\sqrt{hm_j}}b_j.
\]
Since $e^{-ih(\beta_j+2\pi n/h)}=\zeta_j^{-1}$,
we have $\phi_{j,n}\in\Dom(D_F)$ and
\[
 D_F\phi_{j,n}=(\beta_j+2\pi n/h)\phi_{j,n}.
\]
By Lemma \ref{lem_904} and the Fourier basis of $L^2(0,h)$,
the functions $\phi_{j,n}$ form an orthonormal basis of
$L^2((0,h),\mathcal E)$.
Thus $\sigma(D_F)=\Gamma$, and every eigenvalue is simple. 

By the choice of $\beta_j$ in Section \ref{section_2},
the normalized reciprocal root corresponding to $\beta_j$
is $\zeta_j^{-1}$.
The functional equation implies that the multiset of
normalized reciprocal roots is invariant under inversion.
Hence the eigenvalues $\zeta_j$ of multiplication by $F$,
counted with multiplicities $m_j$, agree with the
normalized Frobenius eigenvalues.
The multiplicities $m_j$ are reflected in the
Rosati inner product and 
in the identity $J_\ast(\gamma,\gamma)=1/m_\gamma$ at the endpoint.
The finite-dimensional space $\mathcal E$ is not the
$2g$-dimensional Tate module. It records each distinct
Frobenius eigenvalue once.
Although $M_F$ has eigenspaces of infinite dimension,
the differentiation operator resolves them into simple
Fourier modes.

In Stone's theorem, $D_F$ is the self-adjoint generator
of the strongly continuous unitary group
\eqref{eq_903}.
It is not merely a fixed branch of a logarithm of
Frobenius; its spectrum comprises all lifts of the
Frobenius phases modulo $2\pi/h$.
The value of the group at time $h$ does not determine
its generator uniquely.
It is the full translation group arising from the Weil
quadratic form that selects the operator in
\eqref{eq_904}.

When a Frobenius eigenvalue has multiplicity greater
than one on the Tate module, the present construction
does not, in general, capture the action of geometric
correspondences commuting 
with Frobenius within the corresponding eigenspace.

\section{Comparison with the Riemann zeta function} \label{section_10}

We summarize the correspondence with \cite{Su26a,Su26b}, which
treat the Riemann zeta function.
The construction of the local differentiation operator, the choice of
Riesz vectors $u_{z,a}$, the boundary form, and the identification of
zeros as eigenvalues are common to the two settings.
In both settings, however, the operator representing 
the Weil quadratic form on the usual $L^2$ space must 
be distinguished from the differentiation
operator on the Hilbert space defined by that form.

\subsection{Spectral properties of the local Weil quadratic form}
\hfill \\[5pt]
\noindent\textbf{Friedrichs extension.}
In \cite[Theorem~1.1]{Su26a}, the Friedrichs extension of the
symmetric operator $B_a$, defined from the kernel associated with the
screw function, was identified with the unbounded self-adjoint operator
$A_a$ representing $Q_W$.
Its counterpart here is Proposition~\ref{prop_302}.
Since $A_a$ is bounded, we have the stronger identity
$\overline{B_a}=A_a$.

\medskip
\noindent\textbf{Rayleigh quotients and approximation by test functions.}
\cite[Corollary~1.2]{Su26a} concerns closure in the form norm and
the variational characterization of the lowest eigenvalue.
In the function field case, \eqref{eq_309} follows from
the density of test functions in $L^2$ and boundedness of the
quadratic form.
Corollary~\ref{cor_404} further constructs test functions
that attain the lowest eigenvalue.

\medskip
\noindent\textbf{Continuity of the lowest eigenvalue.}
Continuity of $\lambda_a$ was proved in
\cite[Theorem~1.3]{Su26a}.
For function fields, the finite matrix representation in
Proposition~\ref{prop_403} gives
$\lambda_a=h\lambda_{\min}(H_{\lceil2a/h\rceil})$.
This generally has jumps, and Example~\ref{ex_405} gives a
specific curve for which this occurs.
Moreover, the spectrum of $A_a$ is finite, and every eigenvalue has
infinite multiplicity.
This differs from the discrete spectrum and compact resolvent in the
Riemann zeta case.

\medskip
\noindent\textbf{Asymptotics on small intervals, simplicity, and parity.}
\cite[Theorem~1.4]{Su26a} states that on sufficiently small
intervals the lowest eigenvalue is positive and simple, has an
asymptotic formula with leading term $\log(1/a)$, and has even
eigenfunctions.
For function fields, \eqref{eq_407} gives an exact formula.
The lowest eigenvalue is constant, and its eigenspace is all of
$L^2(-a,a)$.
Neither simplicity nor evenness of all eigenfunctions for the lowest
eigenvalue holds.
Thus the phenomenon of simplicity itself is absent here.

\subsection{The spectrum and limits of the local differentiation operator}

Proposition~\ref{prop_501} and
Theorem~\ref{thm_502} are the function field
counterparts of the construction in \cite[Theorem~1.5]{Su26a}.
The domain of the closure is exactly $H_0^1(-a,a)$, and the adjoint
has the explicit expression $A_a^{-1}i(A_a f)'$.
Unlike $A_a$ itself, the self-adjoint extensions of the
differentiation operator associated with its quadratic form have
compact resolvent and simple eigenvalues.

The convergence of eigenvalues and normalized limits of the
spectral function $W$ in \cite{Su26b} correspond to
Theorems~\ref{thm_1}--\ref{thm_2}.
Both arguments use compression of inverses and the max--min principle.
For function fields, however, sampling is surjective at the finite
value $a_\ast$.
The reproducing kernels therefore have a limit, and even $W$ itself
has an explicit limit.
If there is a central zero, the phase changes, and the local and
global zero eigenspaces must both be removed before compression.
The usual inverse operators cannot be used without this step.

Finally, both constructions first obtain $\Psf$, in which the
multiplicities of the zeros have been removed, and then pass to
$A/E$ using $E_0=E/G$.
Here $\Psf$ records the locations of the zeros, whereas $m_\gamma$
records their multiplicities through the sampling norm and $E_0$.

\subsection{Translation groups and Frobenius}

For the Riemann zeta function, we discussed
the Stone-theoretic interpretation of the global
differentiation operator as the self-adjoint generator
of translations under RH in \cite[Section~7]{Su26a}.
For curves, we obtain the corresponding interpretation
unconditionally in Proposition~\ref{prop_901}.
Furthermore, using the geometric realization in
Section~\ref{section_8}, we identify translation by
$\log q$ with multiplication by normalized Frobenius.
In Theorem~\ref{thm_903}, we also identify the generator
with a differential operator satisfying a Frobenius-twisted
boundary condition.

The essential point is not merely that we can define
a differential operator with a prescribed unitary boundary
condition, but that we identify it with the endpoint
of the local spectral theory arising from the Weil quadratic form.
Thus, in the function field setting, we obtain a geometric
interpretation of analytic structures that also arise
in the Riemann zeta setting of \cite{Su26a,Su26b}.
This does not, however, establish a corresponding
geometric realization in the latter setting.

\appendix

\section{Numerical experiments} \label{section_app}

We compute the spectral function $W(a,\theta_\ast;z)$ explicitly
from arithmetic data given by the numbers of rational points on
curves over finite fields.
In addition to the approximation by finite-dimensional Gram matrices
used for the Riemann zeta function in \cite{Su26b}, the function
field setting permits us to compute $W$ itself as a finite sum
using \eqref{eq_502}.
We compare these two computations.

Put $h=\log q$ and fix $0<a\le a_\ast$.
For $0\le v<1$, put $K(v):=I_a(hv)$.
The set $K(v)$ changes only at points where $h(v+k)=\pm a$
for some $k\in\Z$.
Use these points together with $0,1$ to form a partition
$0=v_0<v_1<\cdots<v_M=1$, and denote by $K_j$ the common set
$K(v)$ for $v\in(v_{j-1},v_j)$.
If $K_j\ne\emptyset$, put
\[
 H_{K_j}:=(p_{\ell-k})_{k,\ell\in K_j}.
\]
Since $K_j$ consists of consecutive integers, relabeling the
indices identifies $H_{K_j}$ with $H_{\#K_j}$ in
\eqref{eq_405}.
Thus $H_{K_j}$ is positive definite by
Theorem~\ref{thm_401}.
Define
\[
 \begin{split}
 \Phi_j(\lambda)
 &:=(v_j-v_{j-1})
 e^{ih\lambda(v_j+v_{j-1})/2}
 \operatorname{sinc}\!\left(
 \frac{h\lambda(v_j-v_{j-1})}{2}\right),
 \end{split}
\]
where $\operatorname{sinc}(z)=\sin z/z$, with
$\operatorname{sinc}(0)=1$.
Substituting $t=hv$ into \eqref{eq_502} and taking the
Fourier transform gives
\begin{equation} \label{eq_a01}
 J_a(w,z)=
 \sum_{\substack{1\le j\le M\\K_j\ne\emptyset}}
 \Phi_j(z-\bar w)
 \sum_{k,\ell\in K_j}
 \bigl(H_{K_j}^{-1}\bigr)_{k\ell}\,
 e^{ih(zk-\bar w\ell)}.
\end{equation}
%
%
Substituting \eqref{eq_a01} into \eqref{eq_510}
therefore computes $W(a,\theta_\ast;z)$ itself from finitely
many inverses of Toeplitz matrices and elementary functions.
\medskip

For comparison, let $N\ge0$ and define
\[
 \chi_{n,a}(x):=
 \frac{e^{i\pi nx/a}}{\sqrt{2a}}\mathbf1_{[-a,a]}(x),
 \qquad
 \mathcal X_{a,N}:=
 \operatorname{span}\{\chi_{n,a}:|n|\le N\},
\]
and the finite-dimensional Gram matrix
\[
 G_{a,N}:=
 \bigl(Q_C(\chi_{m,a},\chi_{n,a})\bigr)_{m,n=-N}^{N}.
\]
By \eqref{eq_306} and
$Q_C(u,f)=\langle A_au,f\rangle_{L^2}$,
\[
 Q_C(\chi_{m,a},\chi_{n,a})
 =h\sum_{|r|h<2a}p_r
 (\chi_{m,a}*\widetilde{\chi_{n,a}})(rh).
\]
The right-hand side can be evaluated by finitely many integrals
of exponential functions over intervals.
Put
$\varphi_{n,a}(z):=\widehat{\chi_{n,a}}(z)
=\sqrt{2a}\,\operatorname{sinc}(az+n\pi)$.
Define
$u_{z,a}^{(N)}
=\sum_{n=-N}^{N}c_n(z)\chi_{n,a}\in\mathcal X_{a,N}$
by $Q_C(v,u_{z,a}^{(N)})=\widehat v(\bar z)$ for
$v\in\mathcal X_{a,N}$.
With the inner product linear in its first variable, this gives
\begin{equation} \label{eq_a02}
 G_{a,N}\overline{\boldsymbol c(z)}
 =\boldsymbol\varphi_{a,N}(\bar z),
 \qquad
 \begin{aligned}
 \boldsymbol c(z)&=(c_n(z))_{n=-N}^{N},\\
 \boldsymbol\varphi_{a,N}(z)
 &=(\varphi_{n,a}(z))_{n=-N}^{N}.
 \end{aligned}
\end{equation}
Solving this system gives
\[
 J_{a,N}(w,z):=\widehat{u_{\bar w,a}^{(N)}}(z),\qquad
 W_N(a,\theta;z):=(z-i)J_{a,N}(-i,z)
       +e^{i\theta}(z+i)J_{a,N}(i,z).
\]
Theorem~\ref{thm_401} gives equivalence of the $Q_C$-norm
and the $L^2$-norm.
Completeness of the exponential system therefore implies that
$\bigcup_N\mathcal X_{a,N}$ is dense in $\Ha$.
Let $\mathsf P_{a,N}$ be the $Q_C$-orthogonal projection onto
$\mathcal X_{a,N}$.
Uniqueness of the Riesz representation gives
$u_{z,a}^{(N)}=\mathsf P_{a,N}u_{z,a}$.
Strong convergence $\mathsf P_{a,N}\to I$ and norm continuity
of $z\mapsto u_{z,a}$ imply that
$u_{z,a}^{(N)}\to u_{z,a}$ uniformly on compact sets of $z$.
Consequently, $J_{a,N}\to J_a$ locally uniformly on $\C^2$,
and $W_N(a,\theta;\cdot)\to W(a,\theta;\cdot)$ locally
uniformly on $\C$.
\medskip

We consider the following three curves:
\[
\begin{aligned}
 C_1/\mathbb F_2 &: \quad y^2+y=x^3+x+1,\\
 C_2/\mathbb F_2 &: \quad y^2+y=x^5+x^3+1,\\
 C_3/\mathbb F_4 &: \quad y^2+y=x^3+\alpha,
 \quad\alpha^2+\alpha+1=0.
\end{aligned}
\]
Counting rational points gives $N_1=1$, $(N_1,N_2)=(1,5)$,
and $N_1=1$, respectively.
Newton's identities for symmetric functions and the functional
equation yield
\[
\begin{aligned}
 P_{C_1}(T)&=1-2T+2T^2,\\
 P_{C_2}(T)&=1-2T+2T^2-4T^3+4T^4,\\
 P_{C_3}(T)&=(1-2T)^2.
\end{aligned}
\]
For $C_1$ and $C_2$, the pairs $(a_\ast,\theta_\ast)$ are
$(\log2,\pi)$ and $(2\log2,\pi)$, respectively;
for $C_3$, the pair is $(\log2,0)$.

Using Mathematica~11, we computed $W$ directly from
\eqref{eq_a01} and $W_N$ from
\eqref{eq_a02}, for $a/a_\ast=1/2,3/4,1$ and
$N=14,16,18$.
The working precision was 70 digits.
Real zeros were found by bracketing sign changes of the
real-valued function
$e^{-i\theta_\ast/2}W(a,\theta_\ast;x)$
and its approximations.
Table~\ref{tab:ff-num-first} compares the first positive zeros
of $W$ and $W_{18}$.

\begin{table}[htbp]
\caption{First positive zeros of $W(a,\theta_\ast;\cdot)$ and
$W_{18}(a,\theta_\ast;\cdot)$, to nine decimal places.}
\label{tab:ff-num-first}
\centering
\small
\setlength{\tabcolsep}{7pt}
\begin{tabular}{ccrr}
\hline
Curve & $a/a_\ast$ & Zero of $W$ & Zero of $W_{18}$\\
\hline
$C_1$ & $1/2$ & 4.532360142 & 4.582550971\\
      & $3/4$ & 1.581275732 & 1.609118244\\
      & $1$   & 1.133090035 & 1.145935686\\
\hline
$C_2$ & $1/2$ & 1.744837513 & 1.764395228\\
      & $3/4$ & 1.281788050 & 1.293732519\\
      & $1$   & 0.377696678 & 0.386740306\\
\hline
$C_3$ & $1/2$ & 9.064720284 & 9.064720284\\
      & $3/4$ & 6.043146856 & 6.043146856\\
      & $1$   & 4.532360142 & 4.532360142\\
\hline
\end{tabular}
\end{table}

For $C_1$ and $C_2$, the computations show the decrease of the
positive eigenvalues as $a\uparrow a_\ast$ and the approach
of the zeros of $W_N$ to those of $W$ as $N$ increases.
For example, the first positive zero of $W_N$ for $C_1$ at the
endpoint is $1.14963542$, $1.14755287$, and $1.14593569$
for $N=14,16,18$, respectively, whereas the first positive zero
of $W$ is $1.13309004\ldots$.
At the endpoint, Theorem~\ref{thm_602} gives the exact identity
$W(a_\ast,\theta_\ast;z)=2B_\ast(i)\Psf(z)/\pi$.
For all three curves, the zeros computed directly from $W$ agree
with $\Gamma$.

The curve $C_3$ has a double central zero.
Here $h=\log4$ and $a_\ast=h/2$, and \eqref{eq_306}
gives $A_a=2hI$ for $0<a\le a_\ast$.
The Riesz representation and the finite-dimensional Fourier
expansion yield, for every $N\ge0$,
\begin{equation} \label{eq_a03}
 \begin{aligned}
 W(a,0;z)&=\frac{2\cosh a}{h}\sin(az),\\
 W_N(a,0;z)&=\frac{2a\sinh a}{h}
 \left(\sum_{n=-N}^{N}\frac{1}{a^2+\pi^2n^2}\right)\sin(az).
 \end{aligned}
\end{equation}
Thus the zeros of $W_N$ and $W$ coincide even for finite $N$
in this example.
On the other hand, $A_{C_3}(z)=-4\sin^2(hz/2)$ has a double
central zero, whereas the zeros of $W(a,0;\cdot)$ are simple.
The ability to compute $W$ itself from finitely many arithmetic
data and to recover $\Gamma$ at the endpoint is a feature that
distinguishes the function field case from the Riemann zeta case.

\bigskip

\noindent
Department of Mathematics, \\
School of Science, \\
Institute of Science Tokyo \\
2-12-1 Ookayama, Meguro-ku, \\
Tokyo 152-8551, Japan \\[2pt]
Email: {\tt msuzuki@math.sci.isct.ac.jp}

\end{document}